\documentclass[sn-mathphys,Numbered]{sn-jnl}% Math and Physical Sciences Reference Style
\usepackage[latin1]{inputenc}
\usepackage[T1]{fontenc}
\usepackage{lmodern}
\usepackage{graphicx}%
\usepackage{multirow}%
\usepackage{amsthm}%
\usepackage{mathrsfs}%
\usepackage[title]{appendix}%
\usepackage{xcolor}%
\usepackage{textcomp}%
\usepackage{manyfoot}%
\usepackage{booktabs}%
\usepackage{algorithm}%
\usepackage{algorithmicx}%
\usepackage{algpseudocode}%
\usepackage{listings}%
\usepackage{amsmath , amsfonts , amssymb , mathrsfs}
\usepackage{stmaryrd}
\usepackage{dsfont}
\usepackage{mathptmx}
\usepackage{parskip}
\usepackage{enumerate}
\usepackage {shadow}
\usepackage{amsmath , amsfonts , amssymb , mathrsfs}

\usepackage[babel]{csquotes}

\usepackage{graphicx}
\usepackage{tocbibind}
\usepackage{microtype}
\usepackage{xcolor}
\usepackage{color}
\usepackage{fancybox}
\usepackage{varioref}
\usepackage{stmaryrd} 
\usepackage{textcomp}
\usepackage{pgf, tikz} 
\usepackage{setspace,booktabs}
\usepackage{enumitem}
\usepackage{pifont}
\usepackage{fancyvrb}

\usepackage{xpatch}  % pour supprimer le point àla fin de theo
\makeatletter   % pour supprimer le point àla fin de theo
\xpatchcmd{\@thm}{\thm@headpunct{.}}{\thm@headpunct{}}{}{}   % pour supprimer le point àla fin de theo
\DeclareSymbolFont{calletters}{OMS}{cmsy}{m}{n}
\DeclareSymbolFontAlphabet{\mathcal}{calletters}
\usepackage{fancyhdr}
\usepackage{fancyvrb}
\theoremstyle{thmstyleone}%
\newtheorem{thm}{Theorem}[section]

\newtheorem{lem}{Lemma}[section]

\newtheorem {ass}{Assumption}[section]

\labelformat{rmq}{Remark~#1}
\labelformat{prop}{Proposition~#1}
\labelformat{lem}{Lemma~#1}
\labelformat{thm}{Theorem~#1}
\labelformat{cor}{Corollary~#1}
\labelformat{ass}{Assumption~#1}

\newtheorem {prop}{Proposition}[section]
\newenvironment{pr}{{\bf{Proof:~}}}
\theoremstyle{thmstyletwo}%

\theoremstyle{thmstylethree}%

\newcommand{\mtt}{\mathtt{T}}
\newcommand{\mds}{\mathds{1}}

\newcommand{\ep}{\varepsilon}
\newcommand{\mt}{\mathsf{T}}

\newcommand{\mteps}{\mathsf{T}_{\!\ep,S}}

\newcommand{\mtepj}{\mathsf{T}_{\!\ep,J}}
\newcommand{\mtj}{\mathsf{T}_{\!\!_J}}
\newcommand{\mtep}{\mathsf{T}_{\ep}}

\newcommand{\un}{\underset}

\newcommand{\n}{\nonumber}
\newcommand{\eps}{\varepsilon}
\newcommand{\displaystylepace}{\mathrm{D}_\eps}

\newcommand{\iep}{\mathcal{I}_{\varepsilon}}
\newcommand{\rep}{\mathcal{R}_{\varepsilon}}
\newcommand{\aep}{\mathcal{A}_{\varepsilon}}

\newcommand{\T}{\mathbb{T}}

\newcommand{\E}{\mathbb{E}}

\renewcommand{\P}{\mathbb{P}}

\newcommand{\bV}{\big\Vert}
\newcommand{\BV}{\Big\Vert}

\newcommand{\nn}{\noindent}
\newcommand{\fpr}{\begin{flushright}
		$\square$
\end{flushright}}
\newcommand{\fprb}{\begin{flushright}
		$\square$
\end{flushright}}

\newcommand{\nbp}{\nabla_{\!\!\ep}^+}
\newcommand{\nbn}{\nabla_{\!\!\ep}^-}
\newcommand{\hd}{\mathtt{H}_{\varepsilon}}
\newcommand{\f}{\mathbf{f}}
\newcommand{\ts}{\mathsf{T}_{\!\ep}}
\newcommand{\mrh}{\mathrm{H}}
\newcommand{\mcru}{\mathscr{U}}
\newcommand{\mcrv}{\mathscr{V}}
\newcommand{\mcrw}{\mathscr{W}}

\newcommand{\scm}{\mathscr{M}}

\newcommand{\bl}{\big\langle \,}
\newcommand{\br}{\,\big\rangle}

\newcommand{\smcr}{\mathscr{M}^S_{\varepsilon}}
\newcommand{\imcr}{\mathscr{M}^I_{\varepsilon}}
\newcommand{\rmcr}{\mathscr{M}^R_{\varepsilon}}

\newcommand{\ungamep}{{_{\mrh^{1-\gamma,\ep}}}}
\newcommand{\ungam}{{_{\mrh^{1-\gamma}}}}

\newcommand{\mgam}{{_{\mrh^{-\gamma}}}}

\newcommand{\mgamz}{{_{\mrh^{-\gamma_0}}}}
\newcommand{\mgamep}{{_{\mrh^{-\gamma , \ep}}}}
\newcommand{\tor}{\mathbb{T}^1}
\newcommand{\dispace}{\mathrm{D}_\eps}

\begin{document}

\title[A SIR epidemic model on a refining spatial grid \\ II-Central limit theorem]{A SIR epidemic model on a refining spatial grid \\ II-Central limit theorem}


\author[1]{\fnm{Thierry} \sur{Gallou\"et}}\email{thierry.gallouet@univ-amu.fr}
\equalcont{These authors contributed equally to this work.}

\author[1]{\fnm{\'Etienne} \sur{Pardoux}}\email{etienne.pardoux@univ-amu.fr}
\equalcont{These authors contributed equally to this work.}

\author*[2]{\fnm{T\'enan} \sur{Yeo}}\email{yeo.tenan@yahoo.fr}

\affil[1]{\orgdiv{Aix Marseille Universit\'e}, \orgname{CNRS, I2M}, \orgaddress{\street{3 Place V. Hugo}, \city{Marseille}, \postcode{13003}, \country{France}}}

\affil*[2]{\orgdiv{Laboratoire de Math\'ematiques Appliqu\'ees et Informatique}, \orgname{Universit\'e F\'elix Houphou\"et-Boigny}, \city{Abidjan}, \postcode{22 BP 582}, \country{C\^ote d'Ivoire}}

%%==================================%%   22 BP 582 Abidjan 22

\abstract{A stochastic  SIR epidemic model taking into account  the heterogeneity of the spatial environment is constructed. The  deterministic model is given by a partial differential equation and the stochastic one by a space-time jump Markov process. The consistency  of the two models is given by a law of large numbers. In this paper, we study the deviation of the spatial stochastic model from the deterministic model by a functional central limit theorem. The limit is a distribution-valued Ornstein-Uhlenbeck Gaussian process, which is the  mild solution of a stochastic partial differential equation.}



\maketitle

\section{Introduction}\label{sec1}

A stochastic spatial model of epidemic has been described by   N'zi et al. (2021) to study the oubreak of infectious diseases in a bounded domain. Such a model takes into account heterogeneity, spatial connectivity and movement of individuals, which play an important role in the spread of the infectious diseases. It is  based on the compartmental SIR model of Kermack and Mckendrick (1927). Let us summarize the results in  N'zi et al. (2021)  in the case of one dimensional space.  

\nn Consider a deterministic and a stochastic SIR model on a grid $\displaystylepace$  of the torus $\tor= [0,1)$ with migration between neighboring sites (two neighboring sites are at distance $\varepsilon$ apart, $ \varepsilon^{-1}\in \mathbb{N}^*$). Let $S_\eps(t,x_i)$ (resp. $I_\eps(t,x_i)$, resp. $R_\eps(t,x_i)$)  be the proportion of the total population which is both susceptible (resp. infectious, resp. removed) and located at site $x_i$ at time
$t$. The dynamics  of susceptible , infected  and removed  individuals at each site can be expressed as
\begin{equation}\label{eqdet} 
\left\{ 
\begin{aligned} 
\dfrac{d\,S_{\varepsilon}}{dt}(t,x_i) &=  \mu_S\,\Delta_{\varepsilon} S_{\varepsilon}(t,x_i)- \dfrac{\beta(x_i)\, S_{\varepsilon}(t,x_i)I_{\varepsilon}(t,x_i)}{S_{\varepsilon}(t,x_i)+I_{\varepsilon}(t,x_i)+R_{\varepsilon}(t,x_i)} \\
\bigskip 
\dfrac{d\,I_{\varepsilon}}{dt}(t,x_i) &= \mu_I\,\Delta_{\varepsilon} I_{\varepsilon}(t,x_i) + \dfrac{\beta(x_i)\, S_{\varepsilon}(t,x_i)I_{\varepsilon}(t,x_i)}{S_{\varepsilon}(t,x_i)+I_{\varepsilon}(t,x_i)+R_{\varepsilon}(t,x_i)}-\alpha(x_i) \,I_{\varepsilon}(t,x_i) \\
\dfrac{d\,R_{\varepsilon}}{dt}(t,x_i) & = \mu_R\,\Delta_{\varepsilon} R_{\varepsilon}(t,x_i)+\alpha(x_i) \,I_{\varepsilon}(t,x_i), \;  \;
(t,x_i)  \in (0,T)\times  \displaystylepace \\
& \hspace{-1.5cm}  S_{\varepsilon}(0,x_i), I_{\varepsilon}(0,x_i), R_{\varepsilon}(0,x_i)\ge 0, \; 0<S_{\varepsilon}(0,x_i)+ I_{\varepsilon}(0,x_i)+ R_{\varepsilon}(0,x_i) \le M,\\
&\hspace{-1.5cm} \text{for some} \; M < \infty ,
\end{aligned}
\right. 
\end{equation}
$\Delta_{\varepsilon}$ is the discrete Laplace operator defined as follows
$$\Delta_{\varepsilon}f(x_i) :=  \varepsilon^{-2}\big[f(x_i+\varepsilon)-2f(x_i)+f(x_i-\varepsilon) \big].$$
The rates $ \beta: [0,1] \longrightarrow  \mathbb{R}_+$ and $ \alpha:[0,1] \longrightarrow  \mathbb{R}_+$ are continuous periodic functions, and $ \mu_S$, $\mu_I$  and  $\mu_R$ are positive diffusion coefficients for the susceptible, infectious and removed subpopulations, respectively. \\
In what follows, we use the notations
$ S_{\varepsilon}(t) := \left(\begin{array}{cl}
S_{\varepsilon}(t,x_1)\\
\vdots \\
S_{\varepsilon}(t,x_{\ell})
\end{array}
\right)$,
$ I_{\varepsilon}(t) := \left( \begin{array}{cl}
I_{\varepsilon}(t,x_1)\\
\vdots \\
I_{\varepsilon}(t,x_{\ell}) 
\end{array} 
\right)$, 
$ R_{\varepsilon}(t) := \left( \begin{array}{cl}
R_{\varepsilon}(t,x_1)\\
\vdots \\
R_{\varepsilon}(t,x_{\ell})
\end{array}
\right)$, and $\displaystyle \quad Z_{\varepsilon}(t)= \big(
S_{\varepsilon}(t)\, ,\,
I_{\varepsilon}(t)\, ,\,
R_{\varepsilon}(t)
\big)^{\mtt}$. Here $\ell = \ep^{-1}$. \\
Note that (\ref{eqdet}) is the discrete space approximation  of the following system of PDEs
\begin{equation}\label{cdm}
\left \{
\begin{aligned}
\dfrac{\partial\,\mathbf{s}}{\partial t}(t,x)= & \mu_S \,\Delta \mathbf{s}(t,x)-\dfrac{\beta(x)\, \mathbf{s}(t,x)\mathbf{i}(t,x) }{\mathbf{s}(t,x)+\mathbf{i}(t,x)+\mathbf{r}(t,x)}\\
\dfrac{\partial\,\mathbf{i}}{\partial t}(t,x)=&\mu_I\, \Delta \mathbf{i}(t,x)+\dfrac{\beta(x)\, \mathbf{s}(t,x)\mathbf{i}(t,x) }{\mathbf{s}(t,x)+\mathbf{i}(t,x)+\mathbf{r}(t,x)} - \alpha(x)\, \mathbf{i}(t,x)\\
\dfrac{\partial\,\mathbf{r}}{\partial t}(t,x)=&\mu_R\, \Delta \mathbf{r}(t,x)+\alpha(x)\, \mathbf{i}(t,x),  \quad (t,x) \in (0,T)\times  D \\
&\hspace{-1.5cm}  \mathbf{s}(0,x), \mathbf{i}(0,x), \mathbf{r}(0,x)\ge 0 , \; 0< \mathbf{s}(0,x)+\mathbf{i}(0,x)+\mathbf{r}(0,x)\le M,
\end{aligned}
\right.
\end{equation}
where $ \displaystyle \Delta = \dfrac{\partial^2}{\partial x^2}$. 
In the sequel, we set  $\mathbf{X}:=(\mathbf{s} \, ,\, \mathbf{i} \, ,\, \mathbf{r})^{\mtt}$.\\
Let $\mathbf{N}$  be the total population size.
The  stochastic version of (\ref{eqdet}) is given by the following system 
\begin{equation}\label{msto}
\hspace{-0.3cm}\left\{ 
\begin{aligned} 
\bigskip 
S_{\mathbf{N},\varepsilon}(t,x_i) &=  S_{\mathbf{N},\varepsilon}(0,x_i) - \frac{1}{\mathbf{N}}\mathrm{P}_{x_i}^{inf}\left( \mathbf{N}\int_0^t \dfrac{\beta(x_i) S_{\mathbf{N},\varepsilon}(r,x_i)I_{\mathbf{N},\varepsilon}(r,x_i)}{S_{\mathbf{N},\varepsilon}(r,x_i)+I_{\mathbf{N},\varepsilon}(r,x_i)+R_{\mathbf{N},\varepsilon}(r,x_i)}dr \right) \\
& \hspace{-1.5cm}- \sum_{y_i\sim x_i}\frac{1}{\mathbf{N}}\mathrm{P}_{S,x_i,y_i}^{mig}\left( \mathbf{N}\int_0^t \frac{\mu_S }{\varepsilon^2}S_{\mathbf{N},\varepsilon}(r,x_i)dr \right) + \sum_{y_i\sim x_i}\frac{1}{\mathbf{N}}\mathrm{P}_{S,y_i,x_i}^{mig}\left(\mathbf{N} \int_0^t \frac{\mu_S }{\varepsilon^2}S_{\mathbf{N},\varepsilon}(r,y_i)dr \right)
\\[2mm] 
I_{\mathbf{N},\varepsilon}(t,x_i) & =  I_{\mathbf{N},\varepsilon}(0,x_i) + \frac{1}{\mathbf{N}}\mathrm{P}_{x_i}^{inf}\left(\mathbf{N}\int_0^t \dfrac{\beta(x_i) S_{\mathbf{N},\varepsilon}(r,x_i)I_{\mathbf{N},\varepsilon}(r,x_i)}{S_{\mathbf{N},\varepsilon}(r,x_i)+I_{\mathbf{N},\varepsilon}(r,x_i)+R_{\mathbf{N},\varepsilon}(r,x_i)}dr \right)  \\
&- \frac{1}{\mathbf{N}}\mathrm{P}_{x_i}^{rec}\left(\mathbf{N}\int_0^t \alpha(x_i)I_{\mathbf{N},\varepsilon}(r,x_i)dr \right) - \sum_{y_i\sim x_i}\frac{1}{\mathbf{N}}\mathrm{P}_{I,x_i,y_i}^{mig}\left(\mathbf{N}\int_0^t\frac{\mu_I }{\varepsilon^2}I_{\mathbf{N},\varepsilon}(r,x_i)dr \right) \\
&+ \sum_{y_i\sim x_i}\frac{1}{\mathbf{N}}\mathrm{P}_{I,y_i,x_i}^{mig}  \left(\mathbf{N}\int_0^t\frac{\mu_I }{\varepsilon^2}I_{\mathbf{N},\varepsilon}(r,y_i)dr \right) 
\\[2mm] 
R_{\mathbf{N},\varepsilon}(t,x_i) & =  R_{\mathbf{N},\varepsilon}(0,x_i) + \frac{1}{\mathbf{N}}\mathrm{P}_{x_i}^{rec}\left(\mathbf{N}\int_0^t \alpha(x_i)I_{\mathbf{N},\varepsilon}(r,x_i)dr \right)\\
&\hspace{-1.5cm}-\sum_{y_i\sim x_i}\frac{1}{\mathbf{N}}\mathrm{P}_{R,x_i,y_i}^{mig}\left(\mathbf{N}\int_0^t \frac{\mu_R}{\varepsilon^2}R_{\mathbf{N},\varepsilon}(r,x_i)dr \right) + \frac{1}{\mathbf{N}}\sum_{y_i\sim x_i}\mathrm{P}_{R,y_i,x_i}^{mig}\left(\mathbf{N}\int_0^t \frac{\mu_R  }{\varepsilon^2}R_{\mathbf{N},\varepsilon}(r,y_i)dr \right)\\
&(t,x_i) \in [0,T]\times \displaystylepace 
\end{aligned}
\right. 
\end{equation}
where all the $\mathrm{P}_j$'s are standard Poisson processes, which are  mutually independent.
For each given site, these  processes count the number of new infectious, recoveries and the migrations between sites.  $ y_i \sim x_i $ means that $ y_i \in \{x_i+\ep \, , \,  x_i-\ep\}$.\\ Let 
$ S_{\mathbf{N},\varepsilon}(t) := \left( \begin{array}{cl}
S_{\mathbf{N},\varepsilon}(t,x_1)\\
\vdots \\
S_{\mathbf{N},\varepsilon}(t,x_{\ell})
\end{array}
\right)$ , 
$ I_{\mathbf{N},\varepsilon}(t) := \left( \begin{array}{cl}
I_{\mathbf{N},\varepsilon}(t,x_1)\\
\vdots \\
I_{\mathbf{N},\varepsilon}(t,x_{\ell}) 
\end{array} 
\right)$, 
$ R_{\mathbf{N},\varepsilon}(t) := \left( \begin{array}{cl}
R_{\mathbf{N},\varepsilon}(t,x_1)\\
\vdots \\
R_{\mathbf{N},\varepsilon}(t,x_{\ell})
\end{array}
\right)$, \\ $\displaystyle \; Z_{\mathbf{N},\varepsilon}(t):= \big(
S_{\mathbf{N},\varepsilon}(t),\,
I_{\mathbf{N},\varepsilon}(t), \,
R_{\mathbf{N},\varepsilon}(t)
\big)^{\mtt}$ and  
$\displaystyle b_{\varepsilon}\big(t,Z_{\mathbf{N},\varepsilon}(t)\big) := \sum_{j=1}^{K}h_j \beta_j( Z_{\mathbf{N}, \varepsilon}(t))$ ($K$ being the number of Poisson's processes in the system), where the vectors  $h_j \in \{ -1 , 0 , 1\}^{3\ell}$ denote the respective jump directions with jump rates $\beta_j$. The SDE (\ref{msto}) can be rewritten as follows
\begin{eqnarray}\label{fsto}
Z_{\mathbf{N},\varepsilon}(t) = Z_{\mathbf{N},\varepsilon}(0)+\int_0^t b_{\varepsilon}\big(r,Z_{\mathbf{N},\varepsilon}(r\big)dr +\dfrac{1}{\mathbf{N}}\sum_{j=1}^{K}h_j \mathrm{P}_j\left(\mathbf{N}\int_{0}^{t}\beta_j\left(Z_{\mathbf{N},\varepsilon}(r)\right)dr\right)  .
\end{eqnarray}
Also the sytem (\ref{eqdet}) can be written as follows
\begin{eqnarray}\label{fdst}
\dfrac{dZ_{\varepsilon}(t)}{dt}= b_{\varepsilon}(t,Z_{\varepsilon}(t)).
\end{eqnarray}

The authors show  the consistency of the two models by a law of large numbers. More precisely,
the  following two results were proved in  N'zi et al. (2021).

\begin{thm}[\textbf{Law  of  Large  Numbers: \boldmath $\mathbf{N}\to \infty$,  $\ep$  being fixed}]\label{llns}~\\
	Let  $ Z_{\mathbf{N},\varepsilon} $ denote the  solution  (\ref{fsto})  and $ Z_{\varepsilon}$ the solution of  (\ref{fdst}). \\
	Let us fix an arbitrary $T > 0$ and assume that $ Z_{\mathbf{N},\varepsilon}(0) \longrightarrow Z_{\varepsilon}(0) $,  as $ \mathbf{N}\rightarrow + \infty $.\\
	Then $ \displaystyle   \underset{0\leq t\leq T}{\sup}\Big\Vert  Z_{\mathbf{N},\varepsilon}(t)-Z_{\varepsilon}(t) \Big\Vert \longrightarrow 0  \; \text{a.s.} \; , \; \; as \; \;   \mathbf{N}\rightarrow + \infty $ . 
\end{thm}
Moreover, for all $ x_i \in D_{\ep}$,  $V_i :=[x_i-\eps/2, x_i+\eps/2 )$ denote the cell centered in the site $x_i$. We define 

$\displaystyle \mathcal{S}_{\varepsilon}(t,x) := \sum_{i=1}^{\varepsilon^{-1}}  S_{\varepsilon}(t,x_i)\mds_{V_i}(x), \; \; \mathcal{I}_{\varepsilon}(t,x) = \sum_{i=1}^{\varepsilon^{-1}} I_{\varepsilon}(t,x_i)\mds_{V_i}(x), \;\;
\mathcal{R}_{\varepsilon}(t,x) := \sum_{i=1}^{\varepsilon^{-1}} R_{\varepsilon}(t,x_i)\mds_{V_i}(x),$

$ \displaystyle \beta(x) := \sum_{i=1}^{\varepsilon^{-1}}\beta(x_i)\mds_{V_i}(x), \;  \alpha(x) := \sum_{i=1}^{\varepsilon^{-1}}\alpha(x_i)\mds_{V_i}(x)$, and we set
\begin{align} \mathbf{X}_{\ep} :=(\mathcal{S}_{\eps}\, ,\, \iep \, ,\, \rep)^{\mtt}. \end{align}

We introduce the canonical projection $P_\eps : L^2(\tor) \longrightarrow \mathtt{H}_\eps$ defined by
$$ f \longmapsto P_\eps f(x)= \eps^{-1}\int_{V_i}f(y)dy, \; \text{if}\; \; x\in V_i .$$
Throughout this paper, we assume that the initial condition satisfies
\begin{ass}\label{init} $\mathbf{s}(0,.)$, $\mathbf{i}(0,.)$, $\mathbf{r}(0,.)$  $\in C^1(\tor)$, \; \;   $\forall x\in \tor$,
	$\mathcal{S}_{\varepsilon}(0,x)=P_\eps \mathbf{s}(0,x)$, $\iep(0,x)=P_\eps \mathbf{i}(0,x)$, $\rep(0,x)=P_\eps \mathbf{r}(0,x)$, and $\displaystyle \int_{\tor} \left(\mathbf{s}(0,x)+\mathbf{i}(0,x)+\mathbf{r}(0,x)\right)dx=1. $ 
\end{ass}
\begin{ass}\label{nu} There exists a constant $c>0$ such that $\displaystyle \inf_{x \in \tor }\mathbf{s}(0,x)\ge c$.
\end{ass}

We use the notation $  \Vert f \Vert_{\infty} := \un{x\in [0, 1]}{\sup}\vert f(x) \vert$ to denote the supremun norm of $f$ in $[0,1]$ and
define $ \Big\Vert\big(f , g, h \big)^{\! \!\mtt} \Big\Vert_{\infty} := \big\Vert f \big\Vert_{\infty} + \big\Vert g \big\Vert_{\infty}+\big\Vert h \big\Vert_{\infty}.$ 

\smallskip

We have the

\begin{thm}\label{csup} For all  $T > 0 $, \; 
	$\displaystyle \un{0\le t\le T}{\sup}\Big\Vert \mathbf{X}_{\varepsilon}(t)- \mathbf{X}(t)\Big\Vert_{\infty} \longrightarrow 0 $, as \,  $ \varepsilon \to 0 . $ 
\end{thm}
\medskip
Next, defining $\displaystyle \mathcal{S}_{\mathbf{N},\varepsilon}(t,x) := \sum_{i=1}^{\varepsilon^{-1}}  S_{\mathbf{N},\varepsilon}(t,x_i)\mds_{V_i}(x), \; \; \mathcal{I}_{\mathbf{N},\varepsilon}(t,x) := \sum_{i=1}^{\varepsilon^{-1}} I_{\mathbf{N},\varepsilon}(t,x_i)\mds_{V_i}(x),$\\
$\displaystyle \mathcal{R}_{\mathbf{N},\varepsilon}(t,x) := \sum_{i=1}^{\varepsilon^{-1}} R_{\mathbf{N},\varepsilon}(t,x_i)\mds_{V_i}(x),$
and setting
$\mathbf{X}_{\mathbf{N},\varepsilon}:=\big(\mathcal{S}_{\mathbf{N},\varepsilon}\, ,\,\mathcal{I}_{\mathbf{N},\varepsilon}\, ,\, \mathcal{R}_{\mathbf{N},\varepsilon}\big)^{\mtt}$, the following theorem is proved  in   N'zi et al. (2021).

\begin{thm}
	Let us assume that $(\ep,\mathbf{N})\to (0,\infty)$, in such way that
	\begin{enumerate}
		\item[(i)] $\dfrac{\mathbf{N}}{\log(1/\varepsilon)}\longrightarrow \infty$ as $\mathbf{N} \to \infty $ and $ \varepsilon\to 0 $;
		\item[(ii)] $\Big\Vert \mathbf{X}_{\mathbf{N},\varepsilon}(0)- \mathbf{X}(0) \Big\Vert_{\infty} \longrightarrow  0 $ in probability.
	\end{enumerate}	
	Then for all $ T > 0 $,  $ \un{0\le t\le T}{\sup}\BV \mathbf{X}_{\mathbf{N},\varepsilon}(t)- \mathbf{X}(t) \BV_{\infty} \longrightarrow 0 $ in probability .
\end{thm}

We devote this paper to study the deviation of the stochastic model from the deterministic one as the mesh size of the grid goes to zero. In this work, we focus  our attention to the periodic boundary conditions on the unit interval $[0,1]$,  which we denote by $\tor$. Let us mention that Blount (1993) and Kotelenez (1986) described similar spatial model for chemical reactions. The resulting process has one component and is compared with the corresponding deterministic model. They proved a functional central limit theorem  under some restriction on the respective speeds of convergence of  the initial number  of particles in each cell and the number  of  cells.\\
The rest of this paper is organized as follows.  In section 2, we give some notations and preliminaries which will be useful in the sequel of this paper.
In section 3, we establish a functional central limit theorem, the main result of this paper, by letting the mesh size $\ep$ of the grid go to zero. The fluctuation limit is a distribution valued generalized Ornstein-Uhlenbeck Gaussian process and can be represented as the solution of a linear stochastic partial differential equation, whose driving terms are Gaussian martingales.

\section{Notations and Preliminaries}

In this section, we give some notations and collect some standard facts on the  Sobolev spaces $\mathrm{H}^{\gamma}(\tor)$, $\gamma \in \mathbb{R}$. First of all, let us  describe some of the properties of the (discrete)-Laplace operator. Let $ \mathtt{H}_{\varepsilon} \subset L^2(\tor)$ denote the space of real valued step functions that are constant on each  cell $ V_i$. 
For $f\in\hd $, let us define
\[ \nbp f(x_i) :=\frac{f(x_i +\ep)-f(x_i)}{\varepsilon}\; \; \text{and}\; \; \nbn f(x_i) :=\frac{f(x_i)- f(x_i-\ep)}{\varepsilon}.  \]
For $ f, g \in L^2(\tor)$, $\displaystyle \langle \; f , g \;  \rangle := \int_{\tor} f(x)g(x)dx $  denotes the scalar product in $\mathrm{L}^2(\tor).$ \\
It is not hard to see that
$$ \langle \; \nbp f, g \; \rangle = -\langle \; f  , \nbn g \; \rangle \; \;  \text{and} \; \; \Delta_{\varepsilon}f=\nbn\nbp f= \nbp\nbn f.$$
For $m$ even and $ x\in \mathbb{R} $ we define
$$
\varphi_{m}(x):=\left\{
\begin{array}{ll}
& 1 , \quad  \quad  \quad\quad\quad \;   \mbox{for} \; \;  m = 0 \\
& \sqrt{2}\cos( m\pi x), \; \text{for} \; \; m\ne 0  \; \text{and even},
\end{array}
\right.
$$
$$
\psi_{m}(x):=\left\{
\begin{array}{ll}
& 0 , \quad  \quad  \quad\quad\quad  \mbox{for} \; \;  m = 0 \\
& \sqrt{2}sin(m\pi x), \; \text{for} \; \; m\ne 0  \; \text{and even}.
\end{array}
\right.
$$
$\left\{  1\, , \, \varphi_m \, , \, \psi_m \, , \, \; m=2k \, , \, k\ge 1  \right\}$ is a complete orthonormal system (CONS) of eigenvectors of $ \Delta $ in $L^2(\tor)$ with eigenvalues  $\displaystyle-\lambda_m =-\pi^2 m^2.$ Consequently, the semigroup $\mathsf{T}(t)\!:= \!\exp(\Delta \,t)$ acting on $L^2(\tor)$  generated by $ \Delta $ can be represented as
$$ \mathsf{T}(t)f\!=\langle \, f,  1 \, \rangle+\!\sum_{k\ge1}\exp(-\lambda_{2k} t )\Big[\langle \; f , \varphi_{2k} \;  \rangle \varphi_{2k}+\langle \; f , \psi_{2k}\;\rangle\psi_{2k} \Big] ,\; \; f\in L^2(\tor).$$
Assume that $ \varepsilon^{-1}$ is an odd integer. 
For $ m \in \left\{ 0 ,2 , \cdots , \varepsilon^{-1}-1 \right\}$, we define $\displaystyle \varphi_{m}^{\varepsilon}(x)=\sqrt{2}\cos(\pi m j\ep)$,  if  $x\in V_j$ and $\displaystyle \psi_{m}^{\varepsilon}(x)=\sqrt{2}\sin(\pi m j\ep)$,  if $x\in V_j$. $\{\, \varphi_{m}^{\varepsilon}, \, \psi_{m}^{\varepsilon}, m \,\} $ form an orthonormal basis of $ \hd$ as a subspace of $L^2\left(\tor\right)$. These vectors  are eigenfunctions of $ \Delta_{\varepsilon} $  with the associated eigenvalues
$\displaystyle -\lambda_{m}^{\varepsilon}=-2\varepsilon^{-2}\big( 1 - \cos(m\pi\varepsilon) \big).$
Note that $ \displaystyle \lambda_{m}^{\varepsilon} \longrightarrow \lambda_{m}$, as $\varepsilon \to 0$. Basic computations show that there exists a constant $ c$,  such that for each  $m$ and $\eps$, $\varepsilon^{-2}\big( 1 - \cos(\pi m \varepsilon) \big) > c\, m^2.$ Let us set $ n_\ep=\frac{\ep^{-1}-1}{2}.$  $ \Delta_{\varepsilon} $ generates a contraction  semigroup $ \ts(t) :=\exp (\Delta_{\varepsilon} t) $  whose action on each $f\in  \hd$ is given by 
\begin{eqnarray}\label{rept}
\ts(t)f= \sum_{k=0}^{n_\ep}\exp(-\lambda_{2k}^{\varepsilon} t )\Big[\langle \; f , \varphi_{2k}^{\varepsilon} \;  \rangle \varphi_{2k}^{\varepsilon}+\langle \; f , \psi_{2k}^{\varepsilon}\;\rangle\psi_{2k}^{\varepsilon} \Big].
\end{eqnarray}
Note that both $ \Delta_{\varepsilon}$ and $\ts(t)$ are self-adjoint and that $\mathsf{T}_{\!\ep}(t)\Delta_{\varepsilon}\varphi=\Delta_{\varepsilon}\ts(t)\varphi.$ \\ For any $ J\in\{S, I, R\}$, the semigroup generated by $ \mu_{\!_J}\Delta$ is $\mt(\mu_{\!_J}t)$. In the sequel, we will use the notation $ \mtj(t):=\mt(\mu_{\!_J}t)$ and similarly, in the discrete case, we will use the notation  $\mtepj(t):=\mtep(\mu_{\!_J}t)$. Also, for any $J\in \{S, I, R\}$, we set $  \lambda_{m,J}:=\mu_{\!_J}\lambda_{m}$  and $ \lambda_{m,J}^{\ep}:=\mu_{\!_J}\lambda_{m}^{\varepsilon}. $\\
For  $ \gamma \in \mathbb{R}_+$, we define the Hilbert space $\mrh^{\gamma}(\tor)$ as follows.
$$ \mrh^{\gamma}(\tor):=\big\{\, f \in L^2\left(\tor\right) , \Vert f \Vert_{_{\mrh^{\gamma}}}^2 := \sum_{m  \; \text{even}}\big[\langle \;f , \varphi_m\;\rangle^2+\langle \;f , \psi_m\;\rangle^2\big](1+\lambda_m)^{\gamma} < \infty \,\big\}.$$ 
We shall use the notations  $\mrh^{\gamma} := \mrh^{\gamma}(\tor)$ and  $L^2 := L^2(\tor)$.\\
Note that $ \displaystyle \Vert \varphi \Vert_{_{\mrh^{\gamma}}}= \Vert( \mathbf{I}-\Delta)^{\gamma/2}\varphi\Vert_{_{L^2}} $, where $ \mathbf{I} $ is the identity operator on $L^2\left(\tor\right).$ For any  three-dimensional vector-valued function $ \Phi =(\Phi_1,\Phi_2,\Phi_3)^{\mtt}$, we  use the notation $\Vert  \Phi \Vert_{_{\mrh^{\gamma}}} := \Big(\Vert  \Phi_1 \Vert_{_{\mrh^{\gamma}}}^2 +\Vert  \Phi_2 \Vert_{_{\mrh^{\gamma}}}^2 +\Vert  \Phi_3 \Vert_{_{\mrh^{\gamma}}}^2\Big)^{1/2} $.\\ For  $ \gamma \in \mathbb{R}$, we also define
$$ \Vert f \Vert_{_{\mrh^{\gamma,\varepsilon}}} := \Bigg[\sum_{m  \; \text{even}} \big(\langle \;f , \varphi_m^{\ep}\;\rangle^2+\langle \;f , \psi_m^{\ep}\;\rangle^2\big)(1+\lambda_m^{\ep})^{\gamma}\Bigg]^{1/2}, \; f \in \hd .$$
For  $f$ , $g \in \hd$,   we have
\begin{eqnarray}\label{sq}
\big\vert \langle f , g\rangle\big\vert &\le& \Vert f\Vert_{_{\mrh^{-\gamma,\varepsilon}}}\Vert g\Vert_{_{\mrh^{\gamma,\varepsilon}}}\; , \quad \gamma\ge 0.
\end{eqnarray}
Elementary calculation shows that for $ f\in \hd$, and $\gamma >0$ there exist  positive constants $c_1(\gamma)$ and $c_2(\gamma)$ such that for all $\eps>0$
\begin{eqnarray}\label{equivl}
c_1(\gamma)\Vert f\Vert_{_{\mrh^{-\gamma,\varepsilon}}}\le \Vert f\Vert_{_{\mrh^{-\gamma}}}\le c_2(\gamma)\Vert f\Vert_{_{\mrh^{-\gamma,\varepsilon}}}.
\end{eqnarray}
$ f^{\prime}:=\dfrac{\partial f}{\partial x }$ will denote the derivative of $f$. \\

\nn In the sequel of this paper we may use the same notation for different constants (we use the generic notation $C$ for a positive constant). These constants can depend upon some  parameters of the model, as long as these are independent of  $ \varepsilon $ and $ \mathbf{N}$,  we will not necessarily mention this dependence explicitly. However, we use $C(\gamma, T)$ to denote a constant which depends  on $\gamma$ and $T$ (and possibly on  some unimportant constants). The exact value may change from line to line.

Let us now consider  the deviation of the stochastic model around its determinsitc law of large numbers limit. To this end we introduce the rescaled difference between $Z_{\mathbf{N},\ep}(t) $ and $Z_{\ep}$, namely 
$$ \Psi_{\mathbf{N},\varepsilon}(t) := \left( \begin{array}{cl}
U_{\mathbf{N},\varepsilon}(t)\\
V_{\mathbf{N},\varepsilon}(t)\\
W_{\mathbf{N},\varepsilon}(t)
\end{array}
\right),$$   
where

$ U_{\mathbf{N},\varepsilon}(t) := \left( \begin{array}{cl}
\sqrt{\mathbf{N}}\Big(S_{\mathbf{N},\varepsilon}(t, x_1)- S_\varepsilon (t,x_1)\Big)\\
\vdots \\
\sqrt{\mathbf{N}}\Big(S_{\mathbf{N},\varepsilon}(t, x_{\ell})- S_\varepsilon (t,x_{\ell})\Big)
\end{array}
\right) $ ,  $ V_{\mathbf{N},\varepsilon}(t) := \left( \begin{array}{cl}
\sqrt{\mathbf{N}}\Big(I_{\mathbf{N},\varepsilon}(t, x_1)- I_\varepsilon (t,x_1)\Big)\\
\vdots \\
\sqrt{\mathbf{N}}\Big(I_{\mathbf{N},\varepsilon}(t, x_{\ell})- I_\varepsilon (t,x_{\ell})\Big)
\end{array}
\right) $ 

and 

\vspace{0.5cm}

$\hspace{3cm} W_{\mathbf{N},\varepsilon}(t) := \left( \begin{array}{cl}
\sqrt{\mathbf{N}}\Big(R_{\mathbf{N},\varepsilon}(t, x_1)- R_\varepsilon (t,x_1)\Big)\\
\vdots \\
\sqrt{\mathbf{N}}\Big(R_{\mathbf{N},\varepsilon}(t, x_{\ell})- R_\varepsilon (t,x_{\ell})\Big)
\end{array}
\right).$ 

In the sequel, we denote by $"\Longrightarrow"$ weak convergence.  By fixing the mesh size $\ep$ of the grid and letting  $\mathbf{N}$ go to infinity, we obtain the following theorem.
\begin{thm}[\textbf{Central Limit Theorem \boldmath : $\mathbf{N}\to \infty $, $\varepsilon $ being fixed}]\label{cltn}~\\
	Assume that $\sqrt{\mathbf{N}}\big( Z_{\mathbf{N},\varepsilon}(0)- Z_{\varepsilon}(0)\big) \longrightarrow 0 $, as $ \mathbf{N}\to \infty $.\\ Then, as $ \mathbf{N} \rightarrow + \infty$ , $\big\{\begin{array}{rl}  \Psi_{\mathbf{N},\varepsilon} (t) ,  \; t\geq 0   \end{array} \big\} \Longrightarrow \big\{\begin{array}{rl} \Psi_{\varepsilon}(t) , \; t\geq 0  \end{array}\big\} ,  $ for the topology of locally uniform convergence, where  the limit process
	$ \Psi_{\varepsilon}(t) := \left( \begin{array}{cl}
	U_{\varepsilon}(t)\\
	V_{\varepsilon}(t)\\
	W_{\varepsilon}(t)
	\end{array}
	\right) $ satisfies
	\begin{equation}\label{oup}
	\Psi_{\varepsilon}(t) = \int_0^t \nabla_{\!\!z} b_{\varepsilon}\big(r,Z_{\varepsilon}(r)\big) \Psi_{\varepsilon}(r)dr + \sum_{j=1}^{K}\int_0^t \sqrt{\beta_j\big(r,Z_{\varepsilon}(r)}\big)dB_j(r) , \quad t \geq 0, 
	\end{equation}
	and $ \displaystyle \{ B_1(t), B_2(t), \cdots, B_K(t) \}$ are mutually independent standard Brownian motions. \\
	More precisely, by setting $ A_{\varepsilon}=S_{\varepsilon} +I_{\varepsilon}+R_{\varepsilon}$, for any site $x_i$, the limit $(U_{\varepsilon}, V_{\varepsilon}, W_{\varepsilon})^{\mtt}$ satisfies the following system
	\begin{eqnarray}
	\hspace{-0.5cm}U_{\varepsilon}(t,x_i)
	&  = & \mu_{S} \int_0^t \Delta_{\varepsilon}U_{\varepsilon}(r,x_i)dr -  \int_0^t \beta(x_i) \dfrac{I_{\varepsilon}(r,x_i)\big(I_{\varepsilon}(r,x_i)+R_{\varepsilon}(r,x_i)\big)V_{\varepsilon}(r,x_i)}{A_{\varepsilon}^2(r,x_i)} dr\nonumber  \\
	&- & \int_0^t\beta(x_i) \dfrac{S_{\varepsilon}(r,x_i)\big(S_{\varepsilon}(r,x_i)+R_{\varepsilon}(r,x_i)\big)U_{\varepsilon}(r,x_i)}{A_{\varepsilon}^2(r,x_i)}  dr
	+\int_0^t \sqrt{\beta(x_i)\dfrac{S_{\varepsilon}(r,x_i)I_{\varepsilon}(r,x_i)}{A_{\varepsilon}(r,x_i)}}\; \; dB_{x_i}^{inf}(r)\nonumber  \\
	&- &\sum_{y_i \sim x_i}\int_0^t\sqrt{\dfrac{\mu_{S}}{\varepsilon ^2} S_{\varepsilon}(r,x_i)}\; \; dB_{x_iy_i}^S(r) 
	+\sum_{y_i \sim x_i}\int_0^t\sqrt{\dfrac{\mu_{S}}{\varepsilon ^2} S_{\varepsilon}(r,y_i)}\; \; dB_{y_ix_i}^S(r) \nonumber
	\end{eqnarray}
\begin{eqnarray}
V_{\varepsilon}(t,x_i)
&=& \mu_{I} \int_0^t \Delta_{\varepsilon}V_{\varepsilon}(r,x_i)dr +\int_0^t \beta(x_i) \dfrac{I_{\varepsilon}(r,x_i)\big(I_{\varepsilon}(r,x_i)+R_{\varepsilon}(r,x_i)\big)V_{\varepsilon}(r,x_i)}{A_{\varepsilon}^2(r,x_i)} dr\nonumber  \\
&+&\int_0^t\beta(x_i) \dfrac{S_{\varepsilon}(r,x_i)\big(S_{\varepsilon}(r,x_i)+R_{\varepsilon}(r,x_i)\big)U_{\varepsilon}(r,x_i)}{A_{\varepsilon}^2(r,x_i)}dr-\int_0^t \alpha(x_i)V_{\varepsilon}(r,x_i)dr\nonumber\\
&-&\int_0^t\sqrt{\beta(x_i)\dfrac{S_{\varepsilon}(r,x_i)I_{\varepsilon}(r,x_i)}{A_{\varepsilon}^2(r,x_i)}}\; \; dB_{x_i}^{inf}(r)+\int_0^t\sqrt{\alpha(x_i) I_{\varepsilon}(r,x_i)}\; \; dB_{x_i}^{rec}(r) \nonumber\\
&-& \sum_{y_i \sim x_i}\int_0^t\sqrt{\dfrac{\mu_{I}}{\varepsilon ^2} I_{\varepsilon}(r,x_i)}dB_{x_iy_i}^I(r)+ \sum_{y_i \sim x_i}\int_0^t\sqrt{\dfrac{\mu_{I}}{\varepsilon ^2} I_{\varepsilon}(r,y_i)}\; \; dB_{y_ix_i}^I(r)  \nonumber 
\end{eqnarray}
\begin{eqnarray}
W_{\varepsilon}(t,x_i) 
& =& \mu_{R} \int_0^t \Delta_{\varepsilon}W_{\varepsilon}(r,x_i)dr+ \int_0^t \alpha(x_i) V_{\varepsilon}(r,x_i)dr-\int_0^t\sqrt{\alpha(x_i) I_{\varepsilon}(r,x_i)}\; \; dB_{x_i}^{rec}(r) \nonumber
\\[2mm]
&-& \sum_{y_i \sim x_i}\int_0^t\sqrt{\dfrac{\mu_{R}}{\varepsilon ^2} R_{\varepsilon}(r,x_i)}\; \;dB_{x_iy_i}^R(r)+ \sum_{y_i \sim x_i}\int_0^t\sqrt{\dfrac{\mu_{R}}{\varepsilon ^2} R_{\varepsilon}(r,y_i)}\; \; dB_{y_ix_i}^R(r), \nonumber 
\end{eqnarray} 
where $\{ B_{x_i}^{inf}: x_i\in \dispace \}$ , $\{ B_{x_i}^{rec} : x_i\in \dispace \}$ , $\{ B_{x_i y_i}^S   : y_i \sim x_i\in \dispace \}$ , $\{ B_{x_i y_i}^I   : y_i \sim x_i\in \dispace \}$ and $\{ B_{x_i y_i}^R   : y_i \sim x_i\in \dispace \}$ are  families of independent Brownian motions.
\end{thm}
\ref{cltn} is a special case of Theorem 3.5 of Kurtz (1971) (see also Theorem 2.3.2 in Britton and Pardoux (2019) ). Then, here, we do not give  the proof and refer the reader to those references for a complete proof. $\square$

Let $ \mathbf{X}=(\mathbf{s} \, ,\, \mathbf{i} \, ,\,\mathbf{r})^{\mtt}$ satisfying the system (\ref{cdm}) on $[0, 1]$. Thanks to Proposition 1.1 of Taylor (1991) (chapter 15, section 1) we have the following lemma.

\begin{lem}\label{unis} Let $\gamma\ge 0 $ and assume that the initial data $ \mathbf{X}(0)$ bebong to  $(\mrh^{\gamma})^3$, then the  parabolic  system (\ref{cdm}) has a unique solution  $ \mathbf{X} \in C\big([0,T]; (\mrh^{\gamma})^3\big) $.
\end{lem}

% In what follows, we use the notation $\nabla_{\!\!\ep}^{+} u:= (\nabla_{\!\!\ep}^{+}u_1,\nabla_{\!\!\ep}^{+}u_2,\nabla_{\!\!\ep}^{+}u_3)^{\mtt}$ for any three-dimensional vector-valued function  $u$, and similarly $ \nabla u :=(\nabla u_1, \nabla u_2, \nabla u_3)^{\mtt}$. 

The rest of this section is devoted to the proof of some estimates for the solution of the system of equations \eqref{eqdet}. We first note that $S_\eps(t,x_i)\ge0, I_\eps(t,x_i)\ge0, R_\eps(t,x_i)\ge0$ for all $t\ge0$, $x_i\in D_\eps$ and $\eps>0$. Moreover for any $T>0$, there exists a contant $C_T$ such that 
\begin{equation}\label{inftyestim}
\sup_{0\le t\le T}\left(\|S_\eps(t)\|_\infty\vee \|I_\eps(t)\|_\infty\vee \|R_\eps(t)\|_\infty\right)\le C_T,\  \forall  \eps>0\,.\end{equation}
Indeed we first note that $\|S_\eps(t)\|_\infty\le M$, since $S_\eps$ is upper bounded by the solution of the ODE
\begin{align} \frac{dX_\eps}{dt}(t,x_i)=\mu_S\Delta_\eps X_\eps(t,x_i),\quad X_\eps(0,x_i)=M\,.\end{align}
Next $I_\eps(t,x_i)$  is upper bounded by the solution of the ODE (with $\bar{\beta}:=\sup_x\beta(x)$) 
\[\frac{dY_\eps}{dt}(t,x_i)=\mu_I\Delta_\eps Y_\eps(t,x_i)+\bar{\beta}Y_\eps(t,x_i),\quad Y_\eps(0,x_i)=M\,.\]
The result for $R_\eps$ is now easy.

Let us set $\mathcal{A}_{\ep} :=\mathcal{S}_{\ep}+\mathcal{I}_{\ep}+\mathcal{R}_{\ep}$ . We have the
\begin{lem}\label{lowbd} 
	For any $T>0$, there exists a positive constant $c_T$ such that \[\mathcal{A}_\eps(t,x)\ge c_T,\quad\text{ for any }\eps>0,\ 0\le t\le T,\ x\in\tor\,. \]
\end{lem}
\begin{pr} 
	We consider the  ODE
	\[ \frac{d\mathcal{S}_\eps}{dt}	(t,x)=\mu_S\Delta_\eps \mathcal{S}_\eps(t,x)
	-\frac{\beta(x)\mathcal{S}_\eps(t,x)\mathcal{I}_\eps(t,x)}{	\mathcal{S}_\eps(t,x)+\mathcal{I}_\eps(t,x)+\mathcal{R}_\eps(t,x)}.\]
	Since $\mathcal{S}_\eps(t,x)+\mathcal{R}_\eps(t,x)\ge0$
	and $\mathcal{I}_\eps(t,x)\ge0$, it is plain that 
	\[ 0\le \frac{\beta(x)\mathcal{I}_\eps(t,x)}{	\mathcal{S}_\eps(t,x)+\mathcal{I}_\eps(t,x)+\mathcal{R}_\eps(t,x)}\le \overline{\beta}, \quad  \text{where}\; \;  \overline{\beta}:=\sup_{x\in \tor}\vert \beta(x)\vert .\]
	Define $\overline{\mathcal{S}}_\eps(t,x)=e^{\overline{\beta}t}\mathcal{S}_\eps(t,x)$. We have
	\[ \frac{d\overline{\mathcal{S}}_\eps}{dt}(t,x)=\mu_S\Delta_\eps \overline{\mathcal{S}}_\eps(t,x)+\left(\overline{\beta}-\frac{\beta(x)\mathcal{I}_\eps(t,x)}{\mathcal{S}_\eps(t,x)+\mathcal{I}_\eps(t,x)+\mathcal{R}_\eps(t,x)}\right)\overline{\mathcal{S}}_\eps(t,x).\]
	Combining this with the last inequality, we deduce that
	\[\overline{\mathcal{S}}_\eps(t,x)\ge [e^{t\mu_S\Delta_\eps} \overline{\mathcal{S}}_\eps(0,\cdot)](x)\ge c,\]
	from  \ref{nu}.
	
	Going back to $\mathcal{S}_\eps$, we note that we have proved that
	\[ \mathcal{S}_\eps(t,x)\ge ce^{-\overline{\beta}t}\,.\]
	In other words, for any $T>0$, there exists a constant $c_T:=ce^{-\overline{\beta}T}$ which is such that
	\[\mathcal{S}_\eps(t,x)\ge c_T,\quad\text{ for any }\eps>0,\ 0\le t\le T,\ x\in\tor\,.\]
	And since $I_\eps(t,x_i)+R_\eps(t,x_i)\ge0$, $\mathcal{A}_\eps(t,x)$ satisfies the same lower bound.
	\fpr
\end{pr}

\begin{lem}\label{estL}
	For any $T>0$, there exists a  constant $C$ such that for each $\eps >0$
	\begin{eqnarray}
	\sup_{0\le t\le T}\Bigg( \bV \mathcal{S}_{\varepsilon}(t)\bV_{L^2}^2+\bV \iep(t)\bV_{L^2}^2+\bV \rep(t)\bV_{L^2}^2\Bigg)\n \\
	&&\hspace{-6cm}+2\int_0^T \Bigg(\mu_S\bV \nabla_{\varepsilon}^+ \mathcal{S}_{\varepsilon}(r) \bV_{L^2}^2 +\mu_I\bV \nabla_{\varepsilon}^+ \iep(r) \bV_{L^2}^2 +\mu_R\bV \nabla_{\varepsilon}^+ \rep(r) \bV_{L^2}^2\Bigg)dr\le C .
	\end{eqnarray}
	
\end{lem}
\begin{pr}  For all $(t,x)\in [0,T]\times [0, 1]  $, we have 
	
	$ \displaystyle \hspace{1cm} \dfrac{d\,\mathcal{S}_{\varepsilon}}{dt}(t,x) =  \mu_S\,\Delta_{\varepsilon} \mathcal{S}_{\varepsilon}(t,x)- \dfrac{\beta(x)\, \mathcal{S}_{\varepsilon}(t,x)\iep(t,x)}{\aep(t,x)}, $ 
	
	which implies
	\begin{eqnarray}
	2\bl \mathcal{S}_{\varepsilon}(t)\,,\, \dfrac{d\,\mathcal{S}_{\varepsilon}}{dt}(t)\br &=&2\mu_S\bl \Delta_{\varepsilon} \mathcal{S}_{\varepsilon}(t)\, , \, \mathcal{S}_{\varepsilon}(t)\br - 2\bl \dfrac{\beta(.)\, \mathcal{S}_{\varepsilon}(t)\iep(t)}{\aep(t)}\, ,\,\mathcal{S}_{\varepsilon}(t) \br\n \\
	&=& -2\mu_S\bl \nabla_{\varepsilon}^+ \mathcal{S}_{\varepsilon}(t)\, , \, \nabla_{\varepsilon}^+ \mathcal{S}_{\varepsilon}(t)\br - 2\bl \dfrac{\beta(.)\, \mathcal{S}_{\varepsilon}(t)\iep(t)}{\aep(t)}\, ,\,\mathcal{S}_{\varepsilon}(t) \br .\n
	\end{eqnarray}
	Then, $\forall\, t \in [0,T]$, 
	\begin{eqnarray}
	\bV \mathcal{S}_{\varepsilon}(t)\bV_{L^2}^2+2\mu_S\int_0^t\bV \nabla_{\varepsilon}^+ \mathcal{S}_{\varepsilon}(r) \bV_{L^2}^2 dr&=& \bV \mathcal{S}_{\varepsilon}(0)\bV_{L^2}^2- 2\int_0^t\bl \dfrac{\beta(.)\, \mathcal{S}_{\varepsilon}(r)\iep(r)}{\aep(r)}\, ,\,\mathcal{S}_{\varepsilon}(r) \br dr.\n
	\end{eqnarray}
	In the same way, we obtain
	\begin{eqnarray}
	\bV \iep(t)\bV_{L^2}^2+2\mu_I\int_0^t\bV \nabla_{\varepsilon}^+ \iep(r) \bV_{L^2}^2 dr&=& \bV \iep(0)\bV_{L^2}^2+ 2\int_0^t\bl \dfrac{\beta(.)\, \mathcal{S}_{\varepsilon}(r)\iep(r)}{\aep(r)}\, ,\,\iep(r) \br dr\n \\
	&&- 2\int_0^t\bl \alpha(.)\iep(r)\, , \, \iep(r) \br dr\, , \n 
	\end{eqnarray}
	and
	\begin{eqnarray}
	\bV \rep(t)\bV_{L^2}^2+2\mu_R\int_0^t\bV \nabla_{\varepsilon}^+ \rep(r) \bV_{L^2}^2 dr&=& \bV \rep(0)\bV_{L^2}^2+ 2\int_0^t\bl \alpha(.)\iep(r)\, , \, \rep(r) \br dr. \n 
	\end{eqnarray}
	Then, we deduce that
	\begin{eqnarray}
	\bV \mathcal{S}_{\varepsilon}(t)\bV_{L^2}^2+\bV \iep(t)\bV_{L^2}^2+\bV \rep(t)\bV_{L^2}^2+2\int_0^t \Bigg(\mu_S\bV \nabla_{\varepsilon}^+ \mathcal{S}_{\varepsilon}(r) \bV_{L^2}^2 +\mu_I\bV \nabla_{\varepsilon}^+ \iep(r) \bV_{L^2}^2 +\mu_R\bV \nabla_{\varepsilon}^+ \rep(r) \bV_{L^2}^2\Bigg)dr\n\\
	&&\hspace{-13cm} \le \bV \mathcal{S}_{\varepsilon}(0)\bV_{L^2}^2+\bV \iep(0)\bV_{L^2}^2+\bV \rep(0)\bV_{L^2}^2+ \int_0^t\Bigg((2\overline{\beta}+\overline{\alpha}) \bV \iep(r)\bV_{L^2}^2+ \overline{\alpha}\bV \rep(r)\bV_{L^2}^2\Bigg) dr, \n
	\end{eqnarray}
	where  $\overline{\alpha}=\un{x\in \tor}{\sup}\vert \alpha(x) \vert$.
	
	It then follows from Gronwall's lemma that
	\begin{eqnarray}
	\bV \mathcal{S}_{\varepsilon}(t)\bV_{L^2}^2+\bV \iep(t)\bV_{L^2}^2+\bV \rep(t)\bV_{L^2}^2+2\int_0^t \Bigg(\mu_S\bV \nabla_{\varepsilon}^+ \mathcal{S}_{\varepsilon}(r) \bV_{L^2}^2 +\mu_I\bV \nabla_{\varepsilon}^+ \iep(r) \bV_{L^2}^2 +\mu_R\bV \nabla_{\varepsilon}^+ \rep(r) \bV_{L^2}^2\Bigg)dr\n\\
	&&\hspace{-9cm} \le \Big(\bV \mathcal{S}_{\varepsilon}(0)\bV_{L^2}^2+\bV \iep(0)\bV_{L^2}^2+\bV \rep(0)\bV_{L^2}^2\Big)e^{C(\overline{\alpha},\overline{\beta})} \n \\
	&&\hspace{-9cm} \le C(\overline{\alpha}, \overline{\beta}) \,.\n
	\end{eqnarray}
	%	Then $\nabla_{\varepsilon}^+ \mathcal{X}$ is bounded in $L^2\big(0,T; (L^2([0,1]))^3\big)$. Hence it converges in $L^2\big(0,T; (L^2([0,1]))^3\big)$ endowed with the weak topology. For any function $\varphi \in \mathcal{C}^1([0,1])$, we have
	%	\begin{eqnarray}
	%	\bl \nabla_{\varepsilon}^+ \mathcal{S}_{\varepsilon} \, , \, \varphi \br&=& -\bl  \mathcal{S}_{\varepsilon} \, , \, \nabla_{\varepsilon}^-\varphi \br\longrightarrow -\bl  \mathbf{s} \, , \, \nabla\varphi \br= \bl  \nabla\mathbf{s} \, , \, \varphi \br
	%	\end{eqnarray}
	%	This relation proves  that $\nabla_{\varepsilon}^+ \mathcal{S}_{\varepsilon} $ converges weakly in $L^2\big(0,T; (L^2([0,1]))^3\big)$ to $\nabla\mathbf{s}$. 
	\fprb
\end{pr}

We now add the following assumption.

\begin{ass}\label{cun}
	The functions $\beta$, $\alpha$ satisfy $\alpha\in C^1(\tor)$ and $\beta\in C^2(\tor)$.
\end{ass} 
Let $f_\ep(t,x):=\beta(x)\frac{\mathcal{S}_{\varepsilon}(t,x)\big[\mathcal{S}_{\varepsilon}(t,x)+\rep(t,x)\big]}{\aep^2(t,x)}, \text{and }  g_\ep(t,x)  := \beta(x)\frac{\iep(t,x)\big[\iep(t,x)+\rep(t,x)\big]}{\aep^2(t,x)}.$

\begin{lem}\label{estgrad} For any $T>0$, there exists a positive constant $C$ such that for all $\eps>0$,
\begin{align} \int_{0}^{T} \left(\bV \nabla_{\varepsilon}^+ f_\ep(t) \bV_{L^2}^2 + \bV \nabla_{\varepsilon}^+ g_\ep(t) \bV_{L^2}^2\right)dt \le C.\end{align}
\end{lem}
\begin{pr}$\forall x\in \tor $, $\forall t\ge 0$ we have
	\begin{eqnarray}\label{sum}
	\nabla_{\varepsilon}^+ f_\ep(t,x)&=& -\dfrac{\beta(x+\ep)\mathcal{S}_{\varepsilon}(t, x+\ep)\big[\mathcal{S}_{\varepsilon}(t, x+\ep)+\rep(t, x+\ep)\big]\big[\aep(t, x+\ep)+\aep(t, x)\big]\nabla_{\ep}^+\aep(t, x)}{\aep^2(t,x)\aep^2(t, x+\ep)}\nonumber\\[2mm]
	&+&\dfrac{\beta(x+\ep)\mathcal{S}_{\varepsilon}(t, x+\ep)}{\aep^2(t, x)}\nabla_{\ep}^+\big(\mathcal{S}_{\varepsilon}(t, x)+\rep(t, x)\big)\nonumber\\[2mm]
	&+&\dfrac{\beta(x+\ep)\big[\mathcal{S}_{\varepsilon}(t, x)+\rep(t, x)\big]}{\aep^2(t, x)}\nabla_{\ep}^+\mathcal{S}_{\varepsilon}(t, x)  \\[2mm]
	&+& \dfrac{\mathcal{S}_{\varepsilon}(t, x)\big[\mathcal{S}_{\varepsilon}(t, x)+\rep(t, x)\big]}{\aep^2(t, x)}\nabla_{\ep}^+\beta(x), \nonumber
	\end{eqnarray}
	
	from which we obtain
	
	\begin{eqnarray}
	\int_{0}^{T}\int_{\tor}\Big\vert\nabla_{\varepsilon}^+ f_\ep(t, x)\Big\vert^2 dx dt &\le&  C \int_{0}^{T}\int_{\tor} \left(\big\vert \nabla_{\ep}^+\beta(x)\big\vert^2+\Big\vert\nabla_{\varepsilon}^+\mathcal{S}_{\varepsilon}(t,x)\Big\vert^2\right. \nonumber\\
	&+&\left.\Big\vert\nabla_{\varepsilon}^+\rep(t,x)\Big\vert^2 +\Big\vert\nabla_{\varepsilon}^+\aep(t,x)\Big\vert^2 \right) dxdt , \nonumber
	\end{eqnarray}
	where we have used  \ref{cun}, inequality \eqref{inftyestim} and \ref{lowbd}. The result now follows from
	\ref{estL}.
	% , we have
	%$$\int_{0}^{T} \int_{\tor} \left(\Big\vert\nabla_{\varepsilon}^+\mathcal{S}_{\varepsilon}(t,x)\Big\vert^2+\Big\vert\nabla_{\varepsilon}^+\rep(t,x)\Big\vert^2+\Big\vert\nabla_{\varepsilon}^+\aep(t,x)\Big\vert^2 \right) dxdt \le C ,$$
	%so, $\displaystyle \int_{0}^{T}\int_{\tor}\Big\vert\nabla_{\varepsilon}^+ f_\ep(t, x)\Big\vert^2 dx dt \le C$.\\ In the same way, we obtain  $\displaystyle \int_{0}^{T}\int_{\tor}\Big\vert\nabla_{\varepsilon}^+ g_\ep(t, x)\Big\vert^2 dx dt \le C$.
	\fpr
\end{pr}

%Let define $\caron{\mathcal{S}}_\ep:=\nabla_{\varepsilon}^+\mathcal{S}_{\varepsilon}$, \;  $\caron{\mathcal{I}}_\ep:=\nabla_{\varepsilon}^+\iep$, \; and $\caron{\mathcal{R}}_\ep:=\nabla_{\varepsilon}^+\rep$.

\begin{lem}\label{hd} For any $T>0$, there exists a positive constant $C$ such that
	\begin{align}
	\sup_{0\le t\le T}\left(\|\nabla^+_\eps S_\eps(t)\|_\infty\vee\|\nabla^+_\eps I_\eps(t)\|_\infty\vee\|\nabla^+_\eps R_\eps(t)\|_\infty\vee\|\nabla^+_\eps f_\eps(t)\|_\infty\vee\|\nabla^+_\eps g_\eps(t)\|_\infty\right)&\le C,\label{estnabinfty}\\
	\int_{0}^{T} \left( \bV \Delta_{\varepsilon} \mathcal{S}_{\varepsilon}(t) \bV_{L^2}^2 +\bV \Delta_{\varepsilon} \iep(t) \bV_{L^2}^2+ \bV \Delta_{\varepsilon} \rep(t) \bV_{L^2}^2\right)dt &\le C,\label{estdelta}\\
	\int_{0}^{T} \left( \bV \Delta_{\varepsilon} f_\ep(t) \bV_{L^2}^2 +\bV \Delta_{\varepsilon} g_{\ep}(t) \bV_{L^2}^2\right)dt &\le C,
	\label{estdeltaf}
	\end{align}
	and 
	\begin{eqnarray}\label{hestfg}
	\sup_{0\le t\le T} \left(\big\Vert f_\ep(t) \big\Vert_{\mrh^{1,\ep}}\vee \big\Vert g_\ep(t) \big\Vert_{\mrh^{1,\ep}}\right)\le C .
	\end{eqnarray}
\end{lem}

%Now, using the estimates of the  above lemmas,  we obtain the
%
%\begin{lem}\label{lfg}
%For any $T>0$, there exists a positive constant $C$ such that $$\int_0^T \left(\Big\Vert \dfrac{d f_\ep}{dt}(t)\Big\Vert_{L^2}^2 + \Big\Vert \dfrac{d g_\ep}{dt}(t)\Big\Vert_{L^2}^2\right)dt \le C.$$
%\end{lem}

\begin{pr} We first etablish \eqref{estnabinfty}. Applying the operator $\nabla^+_\eps$ to the first equation in \eqref{eqdet},  we get
	\begin{align}\label{eqnabla}
	\frac{d\nabla^+_\eps \mathcal{S}_\eps}{dt}(t,x)=\mu_S\Delta_\eps\nabla^+_\eps \mathcal{S}_\eps(t,x)
	-\nabla^+_\eps\left(\frac{\beta\mathcal{S}_\eps \mathcal{I}_\eps}{\mathcal{A}_\eps}\right)(t,x)
	\end{align}
	The last term on the above right hand side is easily explicited thanks to a computation similar to that  done in \eqref{sum}.
	Combining that formula with \ref{cun}, inequality \eqref{inftyestim} and \ref{lowbd}, we deduce that
	\[ \Big\|\nabla^+_\eps\left(\frac{\beta\mathcal{S}_\eps \mathcal{I}_\eps}{\mathcal{A}_\eps}\right)(t)\Big\|_\infty
	\le C\left(\Big\|\nabla^+_\eps \mathcal{S}_\eps(t)\Big\|_\infty+\Big\|\nabla^+_\eps \mathcal{I}_\eps(t)\Big\|_\infty+\Big\|\nabla^+_\eps \mathcal{R}_\eps(t)\Big\|_\infty\right)\,.\]
	From the Duhamel formula,
	\[\nabla^+_\eps \mathcal{S}_\eps(t)=e^{t\mu_S\Delta_\eps}\nabla^+_\eps \mathcal{S}_\eps(0)
	+\int_0^t e^{(t-s)\mu_S\Delta_\eps}\nabla^+_\eps\left(\frac{\beta\mathcal{S}_\eps \mathcal{I}_\eps}{\mathcal{A}_\eps}\right)(s)ds\]
	Since the semigroup $e^{t\mu_S\Delta_\eps}$ is contracting in $L^\infty$, we deduce that
	\[\|\nabla^+_\eps \mathcal{S}_\eps(t)\|_\infty\le \|\nabla^+_\eps \mathcal{S}_\eps(0)\|_\infty+C\int_0^t\left(\|\nabla^+_\eps \mathcal{S}_\eps(s)\|_\infty+\|\nabla^+_\eps \mathcal{I}_\eps(s)\|_\infty+\|\nabla^+_\eps \mathcal{R}_\eps(s)\|_\infty\right)ds\,.\]
	Applying similar arguments to the two other equations in \eqref{eqdet}, we obtain
	\begin{align*}
	\|\nabla^+_\eps \mathcal{S}_\eps(t)\|_\infty
	+\|\nabla^+_\eps \mathcal{I}_\eps(t)\|_\infty
	+\|\nabla^+_\eps \mathcal{R}_\eps(t)\|_\infty
	&\le \|\nabla^+_\eps \mathcal{S}_\eps(0)\|_\infty+\|\nabla^+_\eps \mathcal{I}_\eps(0)\|_\infty
	+\|\nabla^+_\eps \mathcal{R}_\eps(0)\|_\infty\\
	&\ +C\int_0^t\left(\|\nabla^+_\eps \mathcal{S}_\eps(s)\|_\infty+\|\nabla^+_\eps \mathcal{I}_\eps(s)\|_\infty+\|\nabla^+_\eps \mathcal{R}_\eps(s)\|_\infty\right)ds\,.
	\end{align*}
	\eqref{estnabinfty} now follows from Gronwall's Lemma and \ref{init}.
	
	We now multiply \eqref{eqnabla} by $\nabla^+_\eps \mathcal{S}_\eps(t,x)$ and integrate on $[0,t]\times\tor$, yielding
	\begin{align*}
	\|\nabla^+_\eps \mathcal{S}_\eps(t)\|_{L^2}+2\mu_S\int_0^t\|\Delta_\ep\mathcal{S}_\eps(s)\|_{L^2}^2ds
	&=2\int_0^t\left(\frac{\beta\mathcal{S}_\eps \mathcal{I}_\eps}{\mathcal{A}_\eps}(s),\Delta_\eps\mathcal{S}_\eps(s)\right)ds\\
	&\le Ct+\mu_S\int_0^t\|\Delta_\ep\mathcal{S}_\eps(s)\|_{L^2}^2ds\,,
	\end{align*}
	which yields one third of \eqref{estdelta}. The rest of \eqref{estdelta} is proved by similar computations applied to the equations for $\nabla^+_\eps \mathcal{I}_\eps$ and $\nabla^+_\eps \mathcal{R}_\eps$.
	Next \eqref{estdeltaf} follows from \eqref{estdelta}, \eqref{estnabinfty}, \ref{cun}, \eqref{inftyestim} and \ref{lowbd}. 
	
	Since
	$$\big\Vert f_{\ep}(t) \big\Vert_{\mrh^{1,\ep}}^2 \le C\bigg(\big\Vert f_{\ep}(t) \big\Vert_{L^2}^2 +\big\Vert \nabla_{\!\!\ep}^+ f_{\ep}(t)\big\Vert_{L^2}^2\bigg) , $$
	the estimate \eqref{hestfg} follows from \eqref{estnabinfty}, \ref{cun}, inequality \eqref{inftyestim}, \ref{lowbd} and the fact that the norm in $L^2(\tor)$ is bounded by the norm in $L^\infty(\tor)$.
	%Thanks to  the estimates on $\mathcal{S}_{\varepsilon}$, $\iep$, $\rep$ and $\aep$, we have
	%\begin{eqnarray}\label{sld}
	%\sup_{0\le t\le T} \big\Vert f_{\ep}(t) \big\Vert_{L^2}^2 \le C.
	%\end{eqnarray}
	% Moreover, with the previous arguments and exploiting \ref{init} and \ref{cun}, we obtain
	%\begin{eqnarray}\label{sup}
	%\sup_{0\le t\le T}\bigg( \bV \nabla_{\varepsilon}^+\mathcal{S}_{\varepsilon}(t)\bV_{L^2}^2+\bV \nabla_{\varepsilon}^+\iep(t)\bV_{L^2}^2+\bV \nabla_{\varepsilon}^+\rep(t)\bV_{L^2}^2\bigg)
	%&\le& C .
	%\end{eqnarray}
	%Now, it follows from (\ref{sum})  that 
	%\begin{eqnarray}\label{de}
	%\sup_{0\le t\le T}\Big\Vert\nabla_{\varepsilon}^+ f_\ep(t)\Big\Vert_{L^2}&\le& C\sup_{0\le t\le T}\bigg(\Big\Vert \nabla_{\ep}^+\aep(t)\Big\Vert_{L^2}+ \Big\Vert \nabla_{\ep}^+\mathcal{S}_{\varepsilon}(t)\Big\Vert_{L^2}+ \Big\Vert \nabla_{\ep}^+\rep(t)\Big\Vert_{L^2}\bigg) +\Big\Vert \nabla_{\ep}^+\beta\Big\Vert_{L^2}  \nonumber \\
	%&\le& C,
	%\end{eqnarray}
	%where we used (\ref{sup}) for the last inequality.\\
	%Finally, (\ref{sld}) and (\ref{de}) imply that 
	%  $\displaystyle \sup_{0\le t\le T} \big\Vert f_\ep(t) \big\Vert_{\mrh^{1,\ep}}\le C $. 
	%The proof for $ \displaystyle g_\ep $ is similar.
	\fpr
\end{pr}

\begin{lem}\label{convgrad} For any $T>0$, as $\ep \to 0$ 
	$$f_\ep \longrightarrow f , \qquad g_\ep \longrightarrow g, \quad  \nabla_{\ep}^+ f_\ep \longrightarrow \nabla f , \; \;  \text{and} \; \;  \; \; \nabla_{\ep}^+ g_\ep \longrightarrow \nabla g  \; \; \; \; \text{in} \; \; \; \;   C\left([0, T] ; L^2(\T^1)\right),$$
	%$$\Delta_{\ep}f_\ep \longrightarrow \Delta f , \qquad \Delta_{\ep}g_\ep \longrightarrow \Delta g$$ in $L^2\big(]0, T[ \, ; \, L^2\big)$,
	where $\displaystyle f(t,x)=\dfrac{\mathbf{s}(t,x)\left[\mathbf{s}(t,x)+\mathbf{r}(t,x)\right]}{\mathbf{a}^2(t,x)}$\;   and $\;  \displaystyle g(t,x)=\dfrac{\mathbf{i}(t,x)\left[\mathbf{i}(t,x)+\mathbf{r}(t,x)\right]}{\mathbf{a}^2(t,x)}, \; \forall t\in [0, T], \; x\in \tor.$
	
	Moreover $f$, $g \in L^2\left(0, T\; ; \; \mrh^1\right)$.
\end{lem}

\begin{pr} 	
	Let $d$ be the function such that , $\forall$ $ t\in [0, T]$ , $x\in \tor$ and   $\ep>0$ 
	$$ \displaystyle f_\ep(t,x) =d\big(\mathcal{S}_{\varepsilon}(t,x), \iep(t,x), \rep(t,x)\big) \quad \text{and} \quad f(t,x) =d\big(\mathbf{s}(t,x), \mathbf{i}(t,x), \mathbf{r}(t,x)\big).$$
	Furthermore, we know that $\mathcal{S}_{\varepsilon} \longrightarrow \mathbf{s}$, \; $\iep \longrightarrow \mathbf{i}$ and  $\rep \longrightarrow \mathbf{r}$ uniformly on $[0, T] \times \tor$. \\ Since  $d$  is  continuous on $\{(s, i, r)\in (\mathbb{R}_+)^3 : s+i+r>0\}$, then we deduce that $f_\ep \longrightarrow f$ uniformly on $[0, T] \times \tor$, and in particular  in $C\left([0, T] ; L^2(\T^1)\right)$. \\
	From \eqref{eqnabla} and similar equations for  $\displaystyle \nabla_{\!\!\ep}^+ \iep(t,x)$ and $\displaystyle \nabla_{\!\!\ep}^+ \rep(t,x)$, we obtain the convergence of  $\nabla_{\ep}^+ f_\ep \longrightarrow \nabla f$  by an argument similar to the previous one.
	
	The proofs of $g_\ep \longrightarrow g$ and $\nabla_{\ep}^+ g_\ep \longrightarrow \nabla g$ are obtained in the same way.
	%Since $f_\ep$, $g_\ep \in L^2\big(]0, T[ \, ; \, \mrh^{1, \ep}\big)$, and $\displaystyle \dfrac{d f_\ep}{dt}$, $\displaystyle \dfrac{d g_\ep}{dt}\in L^2\big(]0, T[ \, ; \, L^2\big)$, then the proof of the convergence $f_\ep \longrightarrow f$ and  $g_\ep \longrightarrow g$ follows from the Theorem 5.4 in \cite{th}.\\ The proof of the convergence of the other four terms is similar.
	\fpr
\end{pr}

In the sequel, we will write "$f_\eps(t) \longrightarrow f(t) $ in  $\mrh^{1}$"  to mean that "$f_\eps(t) \longrightarrow f(t) $ in  $L^2(\tor)$ and $\nabla_{\!\!\ep}^+f_\eps(t) \longrightarrow \nabla f(t) $ in  $L^2(\tor)$".

%\begin{rmq}
%Note that by using the equation satisfied by $f_\ep$, and using the above estimates, we have  $\dfrac{\partial f_\ep}{\partial t} \longrightarrow \dfrac{\partial f}{\partial t}$ in  $L^2\left(]0, T[, L^2\right)$.
%\end{rmq}

%Let recall the following Moser estimates
%\begin{lem}[Proposition 3.7 of \cite{tay}, page 11]~\\
%	If  $f$ and  $g$ both belong to $ \mrh^{\gamma}$,  with $\gamma>1/2$ then their product also belongs to $\mrh^{\gamma}$ and there exists $C>0$ depending only on $\gamma$ such that $$\big\Vert f g\big\Vert_{_{\mrh^{\gamma}}} \le C \big\Vert f \big\Vert_{_{\mrh^{\gamma}}} \big\Vert g\big\Vert_{_{\mrh^{\gamma}}}.$$
%\end{lem}
We have the  following compactness result.
\begin{lem}[Theorem 1.69  of Bahouri et al. (2011), page 47]~\\
	For any compact subset $E$  of   $\mathbb{R}^d$ and $s_1<s_2$, the embedding of \; $\mrh^{s_2}\left(E\right)$ into $\mrh^{s_1}\left(E\right)$ is a compact linear operator.
\end{lem}

%\begin{lem} ? Let $\gamma <2$.  
%	$f_{\ep}\in L^2\left(0,T ; \mrh^{2}\right)$ , $u_{\ep}\in C\left(0,T ; \mrh^{-\gamma}\right)$, then $f_{\ep}u_{\ep}\in L^2\left(0,T ; \mrh^{-\gamma}\right)$.
%\end{lem}
%\begin{pr}
%	\begin{eqnarray}
%	\int_{0}^{T}\big\Vert f_{\ep}(t)u_{\ep}(t) \big\Vert_{\mrh^{-\gamma}}&=& \int_{0}^{T}\sup_{\substack{\varphi \in \mrh^{\gamma} \\ \Vert \varphi \Vert_{\mrh^{\gamma}}=1}}\langle f_{\ep}(t)u_{\ep}(t) ,\varphi\rangle dt \n\\[2mm]
%	&=& \int_{0}^{T}\sup_{\substack{\varphi \in \mrh^{\gamma} \\ \Vert \varphi \Vert_{\mrh^{\gamma}}=1}}\langle u_{\ep}(t) , f_{\ep}(t)\varphi\rangle \n\\[2mm]
%	&\le&\int_{0}^{T} C\big\Vert f_{\ep}(t)\big\Vert_{\mrh^2}\sup_{\substack{\psi \in \mrh^{\gamma} \\ \Vert \psi \Vert_{\mrh^{\gamma}}=1}}\langle u_{\ep}(t) , \psi\rangle dt\quad \text{(by using the embeding of $\mrh_2$ in $\mrh_{\gamma}$)}\n\\[2mm]
%	&=& C\int_{0}^{T}\big\Vert f_{\ep}(t)\big\Vert_{\mrh^2}\big\Vert u_{\ep}(t)\big\Vert_{\mrh^{-\gamma}}dt\n\\
%	&\le& C_T.\n
%	\end{eqnarray}
%	\fpr
%\end{pr}

In the next section, we  study the behavior of the process $\{\Psi_{\varepsilon}, \; 0<\ep < 1 \}$  as  $\ep$ goes to zero. 

\section{Functional Central Limit Theorem}

Let us define $\displaystyle \mcru_{\varepsilon}(t,x) = \dfrac{1}{\varepsilon^{1/2}}\sum_{i=1}^{\varepsilon^{-1}} U_{\varepsilon}(t,x_i)\mds_{V_i}(x),\; \; $
$\displaystyle\mcrv_{\varepsilon}(t,x) =\dfrac{1}{\varepsilon^{1/2}}\sum_{i=1}^{\varepsilon^{-1}}V_{\varepsilon}(t,x_i)\mds_{V_i}(x),$ 
\[\hspace{0.5cm}\mcrw_{\varepsilon}(t,x) = \dfrac{1}{\varepsilon^{1/2}}\sum_{i=1}^{\varepsilon^{-1}}W_{\varepsilon}(t,x_i)\mds_{V_i}(x).\]
Moreover, we set 
\begin{eqnarray}
\scm_{\varepsilon}^S(t,x) &=& \int_0^t \ep^{-1/2} \sum_{i=1}^{\varepsilon^{-1}}\sqrt{\beta(x_i)\dfrac{S_{\varepsilon}(r,x_i)I_{\varepsilon}(r,x_i)}{A_{\ep}(r,x_i)}}\mds_{V_i}(x)\; \; dB_{x_i}^{inf}(r)\nonumber\\
&+& \sqrt{\mu_S}\int_0^t \ep^{-1/2} \sum_{i=1}^{\varepsilon^{-1}}\sum\limits_{\substack{i\, , \, j \\ x_i \sim  x_j}} \sqrt{ S_{\varepsilon}(r,x_i)}\dfrac{\big( \mds_{V_j}(x)- \mds_{V_i}(x)\big)}{\ep} dB_{x_i x_j}^S(r),\nonumber 
\end{eqnarray}
\begin{eqnarray}
\scm_{\varepsilon}^I(t,x) &=& -\int_0^t \ep^{-1/2} \sum_{i=1}^{\varepsilon^{-1}}\sqrt{\beta(x_i)\dfrac{S_{\varepsilon}(r,x_i)I_{\varepsilon}(r,x_i)}{A_{\ep}(r,x_i)}}\mds_{V_i}(x)\; \; dB_{x_i}^{inf}(r)\nonumber\\
&+&\int_0^t\ep^{-1/2}\sum_{i=1}^{\varepsilon^{-1}}\sqrt{\alpha(x_i) I_{\varepsilon}(r,x_i)}\mds_{V_i}(x)\; dB_{x_i}^{rec}(r) \nonumber\\
&+& \sqrt{\mu_I}\int_0^t \ep^{-1/2} \sum_{i=1}^{\varepsilon^{-1}}\sum\limits_{\substack{i\, , \, j \\ x_i \sim  x_j}} \sqrt{ I_{\varepsilon}(r,x_i)}\dfrac{\big(\mds_{V_j}(x)- \mds_{V_i}(x)\big)}{\ep} dB_{x_i x_j}^S(r),\nonumber
\end{eqnarray}
\begin{eqnarray}
\scm_{\varepsilon}^R(t,x)&=&-\int_0^t\ep^{-1/2}\sum_{i=1}^{\varepsilon^{-1}}\sqrt{\alpha(x_i) I_{\varepsilon}(r,x_i)}\mds_{V_i}(x)\; dB_{x_i}^{rec}(r) \nonumber\\
&+& \sqrt{\mu_R}\int_0^t \ep^{-1/2} \sum_{i=1}^{\varepsilon^{-1}}\sum\limits_{\substack{i\, , \, j \\ x_i \sim  x_j}} \sqrt{ R_{\varepsilon}(r,x_i)}\dfrac{\big(\mds_{V_j}(x)- \mds_{V_i}(x)\big)}{\ep} dB_{x_i x_j}^S(r).\nonumber
\end{eqnarray}
$(\mcru_{\ep}, \mcrv_{\ep}, \mcrw_{\ep})$ satisfies the following system

\begin{equation}\label{syed}
\left\{ 
\begin{aligned} 
\mcru_{\varepsilon}(t)
&=  \int_0^t \mu_S \Delta_{\varepsilon}\mcru_{\varepsilon}(r) dr - \int_0^t \beta(.) \dfrac{\iep(r)\big(\iep(r)+\rep(r)\big)\mcrv_{\varepsilon}(r)}{\aep^2(r)}   dr\\
& \hspace{0.5cm} - \int_0^t  \beta(.) \dfrac{\mathcal{S}_{\varepsilon}(r)\big(\mathcal{S}_{\varepsilon}(r)+\rep(r)\big)\mcru_{\varepsilon}(r)}{\aep^2(r)}   dr + \scm_{\varepsilon}^{S}(t)\\[2mm]
\mcrv_{\varepsilon}(t)
&=  \int_0^t  \mu_{I}\Delta_{\varepsilon} \mcrv_{\varepsilon}(r) dr + \int_0^t  \beta(.) \dfrac{\iep(r)\big(\iep(r)+\rep(r)\big)\mcrv_{\varepsilon}(r)}{\aep^2(r)}   dr\\
& \hspace{-0.5cm} + \int_0^t  \beta(.) \dfrac{\mathcal{S}_{\varepsilon}(r)\big(\mathcal{S}_{\varepsilon}(r)+\rep(r)\big)\mcru_{\varepsilon}(r)}{\aep^2(r)}  dr-\int_0^t \alpha(.)  \mcrv_{\varepsilon}(r) dr+ \scm_{\varepsilon}^{I}(t) \\[2mm]
\mcrw_{\varepsilon}(t)
&=  \int_0^t  \mu_R\Delta_{\varepsilon} \mcrw_{\varepsilon}(r)  dr + \int_0^t \alpha(.) \mcrv_{\varepsilon}(r) dr + \scm_{\varepsilon}^{R}(t). 
\end{aligned}
\right. 
\end{equation}

For $ \gamma \in \mathbb{R}_+$, we denote by $C\big([0,T] ;  \mrh^{-\gamma}\big)$  the complete separable metric space of continuous functions defined on $ [0,T]$  with  values in  $\mrh^{-\gamma}$. For any $ \varepsilon >0$, $\mcru_{\ep}$, $\mcrv_{\varepsilon}$ and $ \mcrw_{\varepsilon}$ can be viewed as continuous processes taking  values in some Hilbert space  $ \mrh^{-\gamma}$. Hence we will study the weak convergence of the process $(\mcru_{\ep}, \mcrv_{\varepsilon}, \mcrw_{\varepsilon})$ in $C\big([0,T] ; (\mrh^{-\gamma})^3\big).$ \\

In the sequel  we will need   to control  the stochastic convolution integrals $\displaystyle  \int_0^t \mtepj(t-r)d\scm_{\varepsilon}^J(r)$, with $J\in \{S, I, R\}$. For that sake, we shall need  a maximal inequality which is a special case of Theorem 2.1 of Kotelenez (1984), which we first recall.

\begin{lem}[Kotelenez (1984)]\label{stopd}  ~\\
	Let $( \mrh \,; \Vert . \Vert_{\mrh}) $ be a separable Hilbert space, $\mathcal{M}$  an $\mrh$-valued locally square integrable càdlàg martingale and  $\mathsf{T}(t)$  a contraction semigroup operator of $ \mathcal{L}(\mrh)$. Then, there 
	is a finite constant $c$ depending only on the Hilbert norm $ \Vert . \Vert_{\mrh} $ such that for all $ T \ge 0 $ 
	\begin{eqnarray}
	\E\bigg(\un{0\le t \le T}{\sup}\Big\Vert \int_0^t \mathsf{T}(t-r)d\mathcal{M}(r)\BV_{\mrh}^2 \bigg)\le c\,e^{4\sigma T}\E\bigg( \BV \mathcal{M}(T)\BV_{\mrh}^2\bigg),
	\end{eqnarray} 
	where $\sigma$ is a real number such that $ \big\Vert\mathsf{T}(t)\big\Vert_{\mathcal{L}(\mrh)} \le e^{\sigma t}$.
\end{lem}

We want to take tke limit as $\epsilon\to 0$ in the system of SDEs (\ref{syed}) satisfied by $\mathscr{Y}_{\ep}$. To this end we will split our system into two subsystems.

First, we consider the following linear system 
\begin{equation}\label{sytlin}
\left\{ 
\begin{aligned} 
d u_{\varepsilon} (t) &= \mu_S\,\Delta_{\ep}u_{\ep}(t) dt+d\smcr(t)\\
dv_{\ep}(t) &=\mu_I\,\Delta_{\varepsilon} v_{\ep}dt+d\imcr(t) \\
dw_{\ep}(t) &=\mu_R\,\Delta_{\varepsilon} w_{\ep}(t)dt+d\rmcr(t)\\
u_{\varepsilon}(0)&=v_{\varepsilon}(0)=w_{\varepsilon}(0)=0.
\end{aligned}
\right. 
\end{equation}
Next, we shall consider the second system
\begin{equation}\label{pb}
\left\{ 
\begin{aligned} 
\dfrac{d\overline{u}_{\varepsilon}}{dt}(t) &=  \mu_S\,\Delta_{\varepsilon} \overline{u}_{\varepsilon}(t)-f_{\varepsilon}(t)\overline{u}_{\varepsilon}(t)-g_{\varepsilon}(t)\overline{v}_{\varepsilon}(t)-f_{\varepsilon}(t)u_{\varepsilon}(t)-g_{\varepsilon}(t)v_{\varepsilon}(t)\\
\dfrac{d\overline{v}_{\varepsilon}}{dt}(t) &=  \mu_I\,\Delta_{\varepsilon} \overline{v}_{\varepsilon}(t)+ f_{\varepsilon}(t)\overline{u}_{\varepsilon}(t)+(g_{\varepsilon}(t)-\alpha)\overline{v}_{\varepsilon}(t)+f_{\varepsilon}(t)u_{\varepsilon}(t)+(g_{\varepsilon}(t)-\alpha)v_{\varepsilon}(t)
\\
\dfrac{d\overline{w}_{\varepsilon}}{dt}(t) &=  \mu_R\,\Delta_{\varepsilon} \overline{w}_{\varepsilon}+ \alpha\left(v_{\varepsilon}+\overline{v}_{\varepsilon}\right)
\\ 
\overline{u}_{\varepsilon}(0)&=\overline{v}_{\varepsilon}(0)=\overline{w}_{\varepsilon}(0)=0,
\end{aligned}
\right. 
\end{equation}
and finally, we note that 
$$\mcru_{\ep}=u_{\ep}+\overline{u}_{\ep}, \quad \mcrv_{\ep}=u_{\ep}+\overline{v}_{\ep}, \quad  \mcrw_{\ep}=w_{\ep}+\overline{w}_{\ep}. $$ 
Then the convergence of  $\mathscr{Y}_\ep :=\left(\mcru_{\ep},  \mcrv_{\ep}, \mcrw_{\ep}\right)$ will follow from both the  convergence of $\left(u_{\ep},  v_{\ep}, w_{\ep}\right)$ and of $\left(\overline{u}_{\ep},  \overline{v}_{\ep}, \overline{w}_{\ep}\right)$.

Let us first look at the convergence of $\left(u_{\ep},  v_{\ep}, w_{\ep}\right)$.

Let $ \scm_{\varepsilon}= \big(\scm_{\varepsilon}^S, \scm_{\varepsilon}^I, \scm_{\varepsilon}^R
\big)^{\mtt}$.

Recall that we denote by $"\Longrightarrow"$ the weak convergence.

\begin{prop}\label{convm} For any $ \gamma> 3/2$, the  Gaussian martingale $ \scm_{\varepsilon}\Longrightarrow \scm := \big(\scm^S, \scm^I,\scm^R\big)^{\mtt}$ in $C\big([0,T]; (\mrh^{-\gamma})^3\big)$ as $\ep \to 0$, where for all $ \varphi \in \mrh^{\gamma}$
	\begin{eqnarray}
	\langle\, \scm^S(t),\varphi \,\rangle \! \! &=&\! \! \! -\int_0^t \int_{\tor} \varphi(x) \sqrt{\dfrac{\beta(x)\mathbf{s}(r,x)\mathbf{i}(r,x)}{\mathbf{a}(r,x)}}\; \dot{W}_1 (dr,dx)  - \sqrt{2\mu_S}\int_0^t\int_{\tor} \varphi^{\prime}(x)\sqrt{\mathbf{s}(r,x)}\; \dot{W}_2(dr,dx), \n 
	\end{eqnarray} 
	\abovedisplayskip = 0.2cm
	\begin{eqnarray}
	\langle\, \scm^I(t),\varphi \,\rangle&=&\int_0^t\int_{\tor}\varphi(x) \sqrt{\dfrac{\beta(x)\mathbf{s}(r,x)\mathbf{i}(r,x)}{\mathbf{a}(r,x)}}\; \dot{W}_1 (dr,dx) +\int_0^t\int_{\tor}\varphi(x)\sqrt{\alpha(x)\mathbf{i}(r,x)}\dot{W}_3(dr,dx) \n \\
	&&-\sqrt{2\mu_I}\int_0^t \int_{\tor}\varphi^{\prime}(x) \sqrt{\mathbf{i}(r,x)}\; \dot{W}_4(dr,dx), \n 
	\end{eqnarray} 
	\abovedisplayskip = 0.1cm 
	\begin{eqnarray}
	\langle\, \scm^R(t),\varphi \,\rangle&=& -\int_0^t\int_{\tor}\varphi(x)\sqrt{\alpha(x)\mathbf{i}(r,x)}\dot{W}_3 (dr,dx)  -\sqrt{2\mu_R}\int_0^t\int_{\tor}\varphi^{\prime}(x)\sqrt{\mathbf{r}(r,x)}\; \dot{W}_5(dr,dx), \n 
	\end{eqnarray} 
	and  $\dot{W_1}$, $\dot{W_2}$, $\dot{W_3}$, $\dot{W_4}$ and $\dot{W_5}$  are standard space-time white noises which are mutually independent.
\end{prop}

\begin{pr} First, we are going to show that there exists a positive constant $C$ independent of $\ep$ such \begin{eqnarray}
	\un{0<\ep<1}{\sup} \E \left(\un{0\le t\le T}{\sup}\big\Vert\scm_{\varepsilon}(t)\big\Vert_{_{\mrh^{-\gamma}}}^2\right) &\le& C.
	\end{eqnarray} 
	Recall that $ \displaystyle \big\Vert\scm_{\varepsilon}(t)\big\Vert_{_{\mrh^{-\gamma}}}^2:= \big\Vert\scm_{\varepsilon}^S(t)\big\Vert_{_{\mrh^{-\gamma}}}^2+\big\Vert\scm_{\varepsilon}^I(t)\big\Vert_{_{\mrh^{-\gamma}}}^2+\big\Vert\scm_{\varepsilon}^R(t)\big\Vert_{_{\mrh^{-\gamma}}}^2$\, . \\ 
	Applying  Doob's  inequality to the martingale $\scm_{\varepsilon}^S$, we have
	\begin{eqnarray}
	\E \left(\un{0\le t\le T}{\sup}\big\Vert \scm_{\varepsilon}^S(t)\big\Vert_{_{\mrh^{-\gamma}}}^2\right)
	&\le& 4 \E \bigg(\big\Vert \scm_{\varepsilon}^S(T)\big\Vert_{_{\mrh^{-\gamma}}}^2\bigg)\n\\
	&& \hspace{-2cm}=4\sum_{m\; \text{even}}\E \Big(\langle \scm_{\varepsilon}^S(T) , \mathbf{f}_{m}\; \rangle^2 \Big)(1+\lambda_m)^{-\gamma}, \; \; \text{with}\; \mathbf{f}_{m}\in \{\varphi_{m}, \psi_{m}\} \nonumber \\
	&&\hspace{-2cm}= \dfrac{4}{\varepsilon}\int_0^T \sum_{m\; \text{even}}\sum_{i=1}^{\varepsilon^{-1}} \dfrac{\beta(x_i)S_{\varepsilon}(r,x_i)I_{\varepsilon}(r,x_i)}{A_{\varepsilon}(r,x_i)}\bigg( \int_{V_i}\mathbf{f}_{m}(x)\,dx \bigg)^2 (1+\lambda_{m})^{-\gamma}dr \nonumber \\
	&&\hspace{-2cm}+\dfrac{4\mu_S}{\varepsilon} \int_0^T \sum_{m \; \text{even}} \sum_{i=1}^{\varepsilon^{-1}}S_{\varepsilon}(r, x_i)\Bigg[\bigg(\int_{V_i} \nabla_{\varepsilon}^{+}\mathbf{f}_{m}(x)\, dx \bigg)^2\!\!+\!\!\bigg(\int_{V_i} \nabla_{\varepsilon}^{-}\mathbf{f}_{m}(x)\, dx \bigg)^2\Bigg](1+\lambda_{m})^{-\gamma} dr. \nonumber
	\end{eqnarray}
	But since $ \dfrac{S_{\varepsilon}(r, x_i)I_{\varepsilon}(r, x_i)}{A_{\varepsilon}(r, x_i)} \le M $ (indeed $\dfrac{I_{\varepsilon}(r, x_i)}{A_{\varepsilon}(r, x_i)}\le 1$ and $S_{\varepsilon}(r, x_i)\le M$, see \eqref{inftyestim} and the line which follows)   \,  and \, $\big\vert \nabla_{\varepsilon}^{\pm}\mathbf{f}_{m}(x) \big\vert^2 \le  2\pi^2 m^2,  $  then we  obtain
	\begin{eqnarray}
	\E \left( \un{0\le t\le T}{\sup}\big\Vert \scm_{\varepsilon}^S(t)\big\Vert_{_{\mrh^{-\gamma}}}^2\right)
	&\le & C(\overline{\beta}, \mu_S, T)\bigg(\sum_{m\; \text{even}}\dfrac{1}{m^{2\gamma}} + \sum_{m\; \text{even}}\dfrac{1}{ m^{2(\gamma-1)}}\bigg).\nonumber
	\end{eqnarray}
	Since $ \displaystyle \sum_{m\; \text{even}}\dfrac{1}{m^{2(\gamma-1)}} < \infty  $ iff $\gamma >3/2$, we then have 
	\begin{eqnarray}\label{majmart}
	\un{0<\ep<1}{\sup} \E \left(\un{0\le t\le T}{\sup}\big\Vert \scm_{\varepsilon}^S(t)\big\Vert_{_{\mrh^{-\gamma}}}^2\right) &\le& C(\overline{\beta}, \mu_S,T), \quad \text{for all $\gamma >3/2$} .
	\end{eqnarray}
	Similar inequalities hold for the martingales $\scm_{\ep}^I$ and $\scm_{\ep}^R$. Hence we obtain
	\begin{eqnarray}\label{mm}
	\un{0<\ep<1}{\sup} \E \left(\un{0\le t\le T}{\sup}\big\Vert \scm_{\varepsilon}(t)\big\Vert_{_{\mrh^{-\gamma}}}^2\right) &\le& C.
	\end{eqnarray}
	Inequality (\ref{mm}) and standard tightness criteria for martingales (see e.g. the proof of \ref{tight}) implies that the  martingale $ \scm_{\ep}$ is tight in $C\big([0,T]; (\mrh^{-\gamma})^3\big)$, with $\gamma>3/2$.\\
	In what follows   $\textbf{<}\textbf{<}\scm_\ep^{S,\gamma_0}\textbf{>}\textbf{>}_t$ denotes the operator--valued increasing process associated to the $L^2(\tor)$--valued martingale $\scm_\ep^{S,\gamma_0}(t)$, whose trace is the increasing process associated to the real valued submartingale $\|\scm_\ep^{S,\gamma_0}(t)\|_{L^2(\tor)}^2$. Let  $ \varphi \in \mrh^{\gamma}$. We set $\scm_{\varepsilon}^{S,\varphi}=\langle \scm_{\varepsilon}^S , \varphi \rangle $. $\scm_{\varepsilon}^{I,\varphi}$ and $\scm_{\varepsilon}^{R,\varphi}$ are defined in the same way.  $\forall \, t\in [0,T]$, we have
	\begin{eqnarray}
	\textbf{<}\textbf{<} \scm_{\ep}^{S,\varphi} \textbf{>}\textbf{>}_t 
	& = & \dfrac{1}{\varepsilon}\int_0^t \sum_{i=1}^{\varepsilon^{-1}} \beta(x_i)\dfrac{S_{\varepsilon}(r,x_i)I_{\varepsilon}(r,x_i)}{A_{\ep}(r,x_i)}\bigg( \int_{V_i}\varphi(x)\,dx \bigg)^2 dr \n \\
	&+& \dfrac{\mu_S}{\varepsilon}\int_0^t \sum_{i=1}^{\varepsilon^{-1}}S_{\varepsilon}(r, x_i)\Bigg[\Big( \int_{V_i} \nabla_{\varepsilon}^{+}\varphi(x)\, dx \Big)^2+\Big( \int_{V_i} \nabla_{\varepsilon}^{-}\varphi(x)\, dx \Big)^2\Bigg] dr. \nonumber
	\end{eqnarray}
	We have 
	\begin{eqnarray}
	\dfrac{1}{\varepsilon}\int_0^t \sum_{i=1}^{\varepsilon^{-1}} \beta(x_i)\dfrac{S_{\varepsilon}(r,x_i)I_{\varepsilon}(r,x_i)}{A_{\varepsilon}(r,x_i)}\bigg( \int_{V_i}\varphi(x)\,dx \bigg)^2 dr
	&&\nonumber \\
	&& \hspace{-6cm} = \dfrac{1}{\varepsilon}\int_0^t \sum_{i=1}^{\varepsilon^{-1}} \beta(x_i)\dfrac{S_{\varepsilon}(r,x_i)I_{\varepsilon}(r,x_i)}{A_{\varepsilon}(r,x_i)}\bigg( \int_{V_i}\varphi(x)\,dx \bigg)\bigg[ \int_{V_i}\Big(\varphi(x)-\varphi(x_i)\Big)\,dx \bigg] dr \nonumber \\
	&& \hspace{-6cm} + \int_0^t \int_{\tor} \sum_{i=1}^{\varepsilon^{-1}} \beta(x_i)\dfrac{S_{\varepsilon}(r,x_i)I_{\varepsilon}(r,x_i)}{A_{\varepsilon}(r,x_i)}\varphi(x)\varphi(x_i)\mds_{V{x_i}}(x) dxdr .\nonumber
	\end{eqnarray}
	On the one hand we have
	\begin{eqnarray}\label{t0}
	\Bigg\vert \dfrac{1}{\varepsilon} \sum_{i=1}^{\varepsilon^{-1}} \beta(x_i)\dfrac{S_{\varepsilon}(r,x_i)I_{\varepsilon}(r,x_i)}{A_{\varepsilon}(r,x_i)}\bigg( \int_{V_i}\varphi(x)\,dx \bigg)\bigg[\int_{V_i}\Big(\varphi(x)-\varphi(x_i)\Big)\,dx \bigg] \Bigg\vert
	&& \nonumber \\
	&& \hspace{-7cm} \leq C\ep \big\Vert \varphi\big\Vert_{\mrh^{\gamma}}\,  \int_{\tor} \dfrac{\beta_{\ep}(x)\mathcal{S}_{\varepsilon}(r,x)\mathcal{I}_{\varepsilon}(r,x)\vert \varphi(x)\vert}{\aep(r,x)} dx \longrightarrow 0,  \n
	\end{eqnarray}
	because the quantity $ \displaystyle \int_{\tor} \dfrac{\beta_{\ep}(x)\mathcal{S}_{\varepsilon}(r,x)\mathcal{I}_{\varepsilon}(r,x)\vert \varphi(x)\vert}{\aep(r,x)} dx$ is bounded uniformly in $ \varepsilon$. \\  Hence \; $\displaystyle  \dfrac{1}{\varepsilon}\int_0^t \sum_{i=1}^{\varepsilon^{-1}} \beta(x_i)\dfrac{S_{\varepsilon}(r,x_i)I_{\varepsilon}(r,x_i)}{A_{\varepsilon}(r,x_i)}\bigg( \int_{V_i}\varphi(x)\,dx \bigg)\bigg[ \int_{V_i}\Big(\varphi(x)-\varphi(x_i)\Big)\,dx \bigg] \longrightarrow 0 $, as $\ep \to 0$. 
	
	On the other hand, the fact that $\displaystyle \un{0\le t\le T}{\sup}\big\Vert \mathbf{X}_{\ep}(t)-\mathbf{X}(t)\big\Vert_{\infty}\longrightarrow 0$, as $\ep \to 0$, leads to 
	\begin{eqnarray}\label{ec} 
	\Bigg\vert  \int_{\tor} \sum_{i=1}^{\varepsilon^{-1}} \beta(x_i)\dfrac{S_{\varepsilon}(r,x_i)I_{\varepsilon}(r,x_i)}{A_{\varepsilon}(r,x_i)}\varphi(x)\varphi(x_i)\mds_{V{x_i}}(x) dx-\int_{\tor} \beta(x) \dfrac{\mathbf{s}(r,x)\mathbf{i}(r,x)}{\mathbf{a}(r,x)}\varphi^2(x) dx\Bigg\vert \longrightarrow 0. \n   
	\end{eqnarray}
	This shows  that 
	$$ \dfrac{1}{\varepsilon}\int_0^t \sum_{i=1}^{\varepsilon^{-1}}\beta(x_i) \dfrac{S_{\varepsilon}(r,x_i)I_{\varepsilon}(r,x_i)}{A_{\varepsilon}(r,x_i)}\bigg( \int_{V_i}\varphi(x)\,dx \bigg)^2 dr\longrightarrow  \int_0^t \int_{\tor} \beta(x) \dfrac{\mathbf{s}(r,x)\mathbf{i}(r,x)}{\mathbf{a}(r,x)}\varphi^2(x) dxdr,$$ as $ \varepsilon\to 0. $ 
	Similar computation shows that 
	$$\dfrac{\mu_S}{\varepsilon}\int_0^t \sum_{i=1}^{\varepsilon^{-1}}S_{\varepsilon}(r, x_i)\Bigg[\Big( \int_{V_i} \nabla_{\varepsilon}^{+}\varphi(x)\, dx \Big)^2+\Big( \int_{V_i} \nabla_{\varepsilon}^{-}\varphi(x)\, dx \Big)^2\Bigg] dr\longrightarrow 2\,\mu_S \int_0^t \int_{\tor} \mathbf{s}(r, x)\big(\varphi^{\prime}(x)\big)^2 dx dr, $$ 
	from which we deduce that
	\begin{eqnarray}\label{cv}
	\textbf{<}\textbf{<} \scm_{\varepsilon}^{S,\varphi} \textbf{>}\textbf{>}_t & \overset{\ep\to 0}{\longrightarrow} &    \int_0^t \int_{\tor} \beta(x) \dfrac{\mathbf{s}(r,x)\mathbf{i}(r,x)}{\mathbf{a}(r,x)}\varphi^2(x) dxdr + 2\,\mu_S \int_0^t \int_{\tor}\mathbf{s}(r,x)\big(\varphi^{\prime}(x)\big)^2dx dr. \nonumber
	\end{eqnarray}
	Hence, if $\dot{W}_1 $ , $\dot{W}_2 $ and $\dot{W}_3 $ are space-time white noises which are mutually independent, so  the limit  of the centered Gaussian martingale $\scm_{\varepsilon}^{S,\varphi}(t) $ can be identified with
	\begin{eqnarray}
	-\int_0^t\int_{\tor}\varphi(x)\sqrt{\dfrac{\beta(x)\mathbf{s}(r,x)\mathbf{i}(r,x)}{\mathbf{a}(r,x)}}\; \dot{W}_1 (dr,dx) -\sqrt{2\mu_S}\int_0^t\int_{\tor}\varphi^{\prime}(x)\sqrt{\mathbf{s}(r,x)}\; \dot{W}_2(dr,dx). \n 
	\end{eqnarray}
	In the same way
	\begin{eqnarray}
	\scm_{\varepsilon}^{I,\varphi}(t)&\!\!\!\Longrightarrow\!\!\!& \int_0^t\int_{\tor}\varphi(x) \sqrt{\dfrac{\beta(x)\mathbf{s}(r,x)\mathbf{i}(r,x)}{\mathbf{a}(r,x)}}\; \dot{W}_1 (dr,dx) +\int_0^t\int_{\tor}\varphi(x)\sqrt{\alpha(x)\mathbf{i}(r,x)}\dot{W}_3(dr,dx) \n \\
	&&-\sqrt{2\mu_I}\int_0^t \int_{\tor}\varphi^{\prime}(x) \sqrt{\mathbf{i}(r,x)}\; \dot{W}_4(dr,dx)\n
	\end{eqnarray}  
	and 
	\begin{eqnarray}
	\scm_{\varepsilon}^{R,\varphi}(t)&\!\!\!\Longrightarrow\!\!\!&  -\int_0^t\int_{\tor}\varphi(x)\sqrt{\alpha(x)\mathbf{i}(r,x)}\dot{W}_3 (dr,dx)-\sqrt{2\mu_R}\int_0^t\int_{\tor}\varphi^{\prime}(x)\sqrt{\mathbf{r}(r,x)}\; \dot{W}_5(dr,dx), \n 
	\end{eqnarray} 
	where  $\dot{W_3}$, $\dot{W_4}$ and $\dot{W_5}$  are also space-time white noises which are mutually independent, and independent from $\dot{W_1}$, $\dot{W_2}$. 
	\fprb	
\end{pr}

Let set $\Im_{\ep}= \big(u_{\ep}\, , \, v_{\ep}\, , \, w_{\ep} \big)^{\mtt}.$\\
We need to check tightness of the sequence of process $\{\Im_{\ep}(t)\, , \, t\in [0, T]\, , 0<\ep < 1 \}$.
\begin{thm}\label{tight}
	For any $\gamma >3/2 $, 
	the  process $\{\Im_{\ep}(t)\, , \, t\in [0, T]\, , 0<\ep < 1 \}$ is tight in   $C\big([0,T]; (\mrh^{-\gamma})^3\big).$ 
\end{thm}
\begin{pr}: We denote by $ \mathcal{G}^T_{\varepsilon}$  the collection of $ \mathcal{F}_t^{\varepsilon}$-stopping times $\overline{\tau}$ such that $\overline{\tau}\le T.$ 
	Following Aldous'  tightness criterion (see Joffe and Metivier (1986)), in oder to show that the process $\{\Im_{\ep}(t)\, , \, t\in [0, T]\, , 0<\ep < 1 \}$ is tight in   $C\big([0,T]; (\mrh^{-\gamma})^3\big)$, it suffices to establish the two following conditions:
	\begin{enumerate}
		\item[\textbf{[T]}] for $\dfrac{3}{2}<\gamma_0< \gamma $,  and $M> 0 $ there exists C such that \;  $\displaystyle \P\Big(\big\Vert \Im_{\ep}(t)\big\Vert_\mgamz\ge M \Big)\le C, $ \; for all $ t\in [0,T]$ ,
		\item[\textbf{[A]}] $\displaystyle
		\un{\theta\to 0}{\lim}\lim\limits_{\substack{\varepsilon\to 0 }}\sup_{\overline{\tau}\in \mathcal{G}_{\varepsilon}^{T-\theta}}\E\bigg(\BV \Im_{\ep}(\overline{\tau}+\theta) - \Im_{\ep}(\overline{\tau})\BV_\mgam^2\bigg) = 0 . $
	\end{enumerate}
	
	Let $\dfrac{3}{2}\! < \! \gamma_0\!< \!\gamma $. 
	Let us set $ \displaystyle u_{\ep}^{\gamma_0}(t,x)=(\mathbf{I}-\Delta_{\ep})^{-\gamma_0/2}u_{\ep}(t,x)$.  $ \forall \, t\in [0,T]$, we have $$ \big\Vert u_{\ep}(t)\big\Vert_\mgamz^2=\big\langle\, u_{\ep}^{\gamma_0}(t), u_{\ep}^{\gamma_0}(t) \,\big\rangle.$$  
	If we define $\scm_\ep^{S,\gamma_0}(t):=(\mathbf{I}-\Delta_{\ep})^{-\gamma_0/2}\scm_\ep^S(t)$, since $\gamma_0>3/2$, it follows from (\ref{convm}) that $\scm_\ep^{S,\gamma_0}(t)$ is bounded as $\ep\to0$, as an $L^2(\tor)$--valued 
	martingale.
	Applying the Itô formula to $\vert u_{\eps}^{\gamma_0}(t,x) \vert^2$ and integrating over $\tor$ leads to 
	\begin{eqnarray}
	\big\Vert u_{\ep}(t)\big\Vert_\mgamz^2&=& -2\int_0^t \big \langle\, \nabla_{\!\!\ep}^{+} u_{\ep}^{\gamma_0}(r), \nabla_{\!\!\ep}^{+} u_{\ep}^{\gamma_0}(r)\, \big\rangle dr
	+ 2\int_0^t\big\langle\, u_{\ep}^{\gamma_0}(r),d\scm_\ep^{S,\gamma_0}(r) \,\big\rangle \n\\
	&+& \int_{\tor}\textbf{<}\textbf{<} \scm_\ep^{S,\gamma_0}(.,x) \textbf{>}\textbf{>}_t dx.\n
	\end{eqnarray}
	Letting $t=T$ and taking the expectation, we deduce that
	\begin{align*}
	\E(\big\Vert u_{\ep}(T)\big\Vert_\mgamz^2)+2\mu_S\E\int_0^T\big\Vert \nabla_{\!\!\ep}^+u_{\ep}(t)\big\Vert_\mgamz^2dt
	=\E\left(\|\scm_\ep^{S,\gamma_0}(T)\|^2_{L^2}\right)\,.
	\end{align*}
	Next we want to take the supremum on $[0, T]$ in the previous identity. For that sake, we use the Burkholder Davis Gundy inequality, which implies that
	\begin{align*}
	\E\left[\sup_{0\le t\le T}\left|\int_0^t\big\langle\, u_{\ep}^{\gamma_0}(r),d\scm_\ep^{S,\gamma_0}(r) \,\big\rangle\right|\right]&\le 3\E\sqrt{\textbf{<}\textbf{<}\int_0^\cdot\big\langle\,u_{\ep}^{\gamma_0}(r),d\scm_\ep^{S,\gamma_0}(r) \,\big\rangle \textbf{>}\textbf{>}_T}\\
	&\le3\E\left(\sup_{0\le t\le T}\|u_{\ep}^{\gamma_0}(t)\|_{L^2}\sqrt{Tr\textbf{<}\textbf{<}\scm_\ep^{S,\gamma_0}\textbf{>}\textbf{>}_T}\right)\\
	&\le \frac{1}{2}\E\left(\sup_{0\le t\le T}\|u_{\ep}^{\gamma_0}(t)\|_{L^2}^2\right)+\frac{9}{2}\E(\|\scm_\ep^{S,\gamma_0}(T)\|_{L^2}^2).
	\end{align*}
	%where $\textbf{<}\textbf{<}\scm_\ep^{S,\gamma_0}\textbf{>}\textbf{>}_t$ denotes the operator--valued increasing process associated to the $L^2(\tor)$--valued martingale $\scm_\ep^{S,\gamma_0}(t)$, whose trace is the increasing process associated to the real valued submartingale $\|\scm_\ep^{S,\gamma_0}(t)\|_{L^2(\tor)}^2$. 
	We then obtain, thanks to \eqref{majmart},
	\begin{align*}
	\E\left(\sup_{0\le t\le T}\big\Vert u_{\ep}(t)\big\Vert_\mgamz^2\right) &= 11\ \E\left(\sup_{0\le t\le T}\big\Vert \scm_\ep^{S,\gamma_0}(t)\big\Vert_{L^2}^2\right)\n \\
	&\le   44\ C(\overline{\beta}, \mu_S,T)\,.
	\end{align*}
	
	We also obtain similar inequalities for $v_\ep$ and $w_\ep$. Hence there exists a constant $C$ such that for all $\ep>0$,
	\begin{align}
	\E&\left(\sup_{0\le t\le T}\big\Vert u_{\ep}(t)\big\Vert_\mgamz^2 +\sup_{0\le t\le T}\big\Vert v_{\ep}(t)\big\Vert_\mgamz^2+\sup_{0\le t\le T}\big\Vert w_{\ep}(t)\big\Vert_\mgamz^2\right)\nonumber\\
	&\quad+ 2\E\int_{0}^{T}\left[\mu_S\big\Vert \nabla_{\!\!\ep}^+u_{\ep}(r)\big\Vert_\mgamz^2 
	+ \mu_I \big\Vert \nabla_{\!\!\ep}^+v_{\ep}(r)\big\Vert_\mgamz^2+\mu_R\big\Vert \nabla_{\!\!\ep}^+w_{\ep}(r)\big\Vert_\mgamz^2\right]dr\le C\,.\label{nabl}
	\end{align}
	Then \textbf{[T]} follows by using Markov's inequality.

	Let $\theta >0$ and $ \overline{\tau} \in \mathcal{G}^{T-\theta}_{\varepsilon} $. We have
	\begin{eqnarray}
	u_{\varepsilon}(\overline{\tau}+\theta)-u_{\varepsilon}(\overline{\tau})&=& \big[ \mteps(\theta)-\mathbf{I}\big]u_{\varepsilon}(\overline{\tau})+ \int_{\overline{\tau}}^{\overline{\tau}+\theta}\mteps(\overline{\tau}+\theta-r)d\scm_{\varepsilon}^S(r). \nonumber
	\end{eqnarray}
	So,
	\begin{eqnarray}
	\mathbb{E}\Bigg(\Big\Vert u_{\varepsilon}(\overline{\tau} + \theta)- u_{\varepsilon}(\overline{\tau})\Big\Vert_\mgam^2\Bigg)
	&\le& 2\E\Bigg(\Big\Vert \big[ \mteps(\theta)-\mathbf{I}\big]u_{\varepsilon}(\overline{\tau}) \Big\Vert_\mgam^2 \Bigg)\n\\
	&+&2\E\Bigg( \Big\Vert \int_{\overline{\tau}}^{\overline{\tau}+\theta} \mteps(\overline{\tau}+\theta-r) d\scm_{\varepsilon}^S(r)  \Big\Vert_\mgam^2\Bigg).\nonumber 
	\end{eqnarray}
	Let us deal with each term separately. 
	First using the inequality (\ref{equivl}), there is a constant $ C(\gamma)$ such that
	\begin{eqnarray}
	\E\bigg(\Big\Vert \big[ \mteps(\theta)-\mathbf{I}\big] u_{\varepsilon}(\overline{\tau})\Big\Vert_\mgam^2\bigg)
	&\le & C(\gamma) \E\bigg(\Big\Vert \big[ \mteps(\theta)-\mathbf{I}\big]u_{\varepsilon}(\overline{\tau})\Big\Vert_\mgamep^2\bigg).\nonumber
	\end{eqnarray} 
	Let $ 3/2<\gamma^{\prime}<\gamma$, and let $c$  a positive constant. We have
	\begin{eqnarray}
	\Big\Vert\big[\mteps(\theta)-\mathbf{I}\big] u_{\varepsilon}(\overline{\tau})\Big\Vert_\mgam^2 &=& \sum_{\lambda_m^\ep\ge c} \bl \big[ \mteps(\theta)-\mathbf{I}\big]u_{\varepsilon}(\overline{\tau}), \f_{m}^{\ep} \br^2\big( 1+\lambda_{m}^{\varepsilon}\big)^{-\gamma} \n \\
	&+& \sum_{\lambda_m^\ep< c} \bl \big[ \mteps(\theta)-\mathbf{I}\big]u_{\varepsilon}(\overline{\tau}), \f_{m}^{\ep} \br^2\big( 1+\lambda_{m}^{\varepsilon}\big)^{-\gamma} ,\n
	\end{eqnarray}
	and
	\begin{eqnarray}
	\sum_{\lambda_m^\ep\ge c} \bl \big[ \mteps(\theta)-\mathbf{I}\big]u_{\varepsilon}(\overline{\tau}), \f_{m}^{\ep} \br^2\big( 1+\lambda_{m}^{\varepsilon}\big)^{-\gamma}\n\\
	&&\hspace{-5cm} \le
	( 1+c)^{\gamma^{\prime}-\gamma}\sum_{\lambda_m^\ep\ge c} \bl \big[\mteps(\theta)-\mathbf{I}\big] u_{\varepsilon}(\overline{\tau}), \f_{m}^{\ep} \br^2\big(1+\lambda_{m}^{\varepsilon}\big)^{-\gamma^{\prime}}\n\\
	&& \hspace{-5cm} \le (1+c)^{\gamma^{\prime}-\gamma} \Big\Vert\big[\mteps(\theta)-\mathbf{I}\big] u_{\varepsilon}(\overline{\tau})\Big\Vert_{\mrh^{-\gamma^{\prime}}}^2.\n
	\end{eqnarray} Then
	\begin{eqnarray}
	\E\bigg(\Big\Vert \big[ \mteps(\theta)-\mathbf{I}\big]u_{\varepsilon}(\overline{\tau})\Big\Vert_\mgam^2\bigg)
	&\le& C(\gamma)(1+c)^{\gamma^{\prime}-\gamma} \E\left(\Big\Vert\big[\mteps(\theta)-\mathbf{I}\big] u_{\varepsilon}(\overline{\tau})\Big\Vert_{\mrh^{-\gamma^{\prime}}}^2\right)\n \\
	&&+ C(\gamma) \sum_{\lambda_m^{\ep}<c} \big(e^{-\lambda_{m}^{\varepsilon}\theta}-1\big)^2\E\Big( \langle \; u_{\varepsilon}(\overline{\tau}), \mathbf{f}_m^{\varepsilon}\; \rangle^2 \Big)\big( 1+\lambda_{m}^{\varepsilon}\big)^{-\gamma}\n
	\end{eqnarray}
	On the one hand, since  $\E\left(\Big\Vert\big[\mteps(\theta)-\mathbf{I}\big] u_{\varepsilon}(\overline{\tau})\Big\Vert_{\mrh^{-\gamma^{\prime}}}^2\right)\le C$, we can choose $c$ large enough such that $\displaystyle C(\gamma)(1+c)^{\gamma^{\prime}-\gamma} \E\left(\Big\Vert\big[\mteps(\theta)-\mathbf{I}\big] u_{\varepsilon}(\overline{\tau})\Big\Vert_{\mrh^{-\gamma^{\prime}}}^2\right)\le \ep/2.$\\
	On the other hand, we have
	\begin{eqnarray}
	\sum_{\lambda_m^{\ep}<c} \big(e^{-\lambda_{m}^{\varepsilon}\theta}-1\big)^2\E\Big( \langle \;u_{\varepsilon}(\overline{\tau}), \mathbf{f}_m^{\varepsilon}\; \rangle^2 \Big)\big( 1+\lambda_{m}^{\varepsilon}\big)^{-\gamma}\n\\
	&&\hspace{-5cm}\le\un{\lambda_m^{\ep}<c}{\sup}\big(1-e^{-\lambda_n^{\varepsilon}\theta}\big)^2\sum_{\lambda_m^{\ep}<c} \E\Big(\langle \;u_{\varepsilon}(\overline{\tau}), \mathbf{f}_m^{\varepsilon}\; \rangle^2 \Big)\big(1+\lambda_{m}^{\varepsilon}\big)^{-\gamma}\nonumber\\
	&&\hspace{-5cm}\le \un{\lambda_m^{\ep}<c}{\sup}\big(1-e^{-\lambda_m^{\varepsilon}\theta}\big)^2\mathbb{E}\bigg(\big\Vert u_{\varepsilon}(\overline{\tau})\big\Vert_\mgamep^2\bigg)\nonumber\\
	&& \hspace{-5cm} \le C(\gamma)\un{\lambda_m^{\ep}<c}{\sup}\big(1-e^{-\lambda_m^{\varepsilon}\theta}\big)^2 \mathbb{E}\bigg(\big\Vert u_{\varepsilon}(\overline{\tau})\big\Vert_\mgam^2\bigg). \nonumber
	\end{eqnarray}
	Since 
	\begin{eqnarray}
	\mathbb{E}\bigg(\big\Vert u_{\varepsilon}(\overline{\tau})\big\Vert_{-\gamma}^2\bigg) &\le& \mathbb{E}\bigg(\sup_{0\le t\le T}\big\Vert u_{\varepsilon}(t)\big\Vert_{-\gamma}^2\bigg) \n \\
	&\le &  C ,\n
	\end{eqnarray}
	then for the previous choice of $c$, we can choose $\theta$ small enough such that $C(\gamma)\un{\lambda_m^{\ep}<c}{\sup}\big(1-e^{-\lambda_m^{\varepsilon}\theta}\big)^2 \mathbb{E}\bigg(\big\Vert u_{\varepsilon}(\overline{\tau})\big\Vert_\mgam^2\bigg)\le \ep/2$. Hence
	\begin{align*}
	\un{\theta\to 0}{\lim}\lim\limits_{\substack{\varepsilon\to 0 }}\sup_{\overline{\tau}\in \mathcal{G}_{\varepsilon}^{T-\theta}}\E\bigg(\Big\Vert \big[ \mteps(\theta)-\mathbf{I}\big]u_{\varepsilon}(\overline{\tau})\Big\Vert_\mgam^2\bigg) = 0\,. 
	\end{align*}

	Secondly, using the equivalence of the norms $\Vert . \Vert_\mgam$ and $\Vert  . \Vert_\mgamep$\, , and the fact that $\mteps$ is a contraction semigroup on $ \mathtt{H}_{\ep}$ we have 
	\begin{eqnarray}\label{bmart}
	\E\Bigg(\Big\Vert \int_{\overline{\tau}}^{\overline{\tau}+\theta} \mteps(\overline{\tau}+\theta-r) \scm_{\varepsilon}^S(r)  \Big\Vert_\mgam^2 \Bigg)
	& = &\E\Bigg( \Big\Vert \int_0^{\theta} \mteps(\theta-r) d\scm_{\varepsilon}^S(r+\overline{\tau})  \Big\Vert_\mgam^2\Bigg) \nonumber \\
	&& \hspace{-4cm}\le C(\gamma)\mathbb{E}\Bigg( \Big\Vert \scm_{\varepsilon}^S(\overline{\tau} + \theta)- \scm_{\varepsilon}^S(\overline{\tau})\Big\Vert_\mgamep^2\Bigg)\nonumber\\
	& & \hspace{-4cm}\le 2 C(\gamma) \E\Bigg(\Big\Vert \int_{\overline{\tau}}^{\overline{\tau}+\theta} \sum_{i=1}^{\varepsilon^{-1}}\dfrac{\sqrt{\beta(x_i)}}{\varepsilon}\sqrt{\dfrac{ S_{\varepsilon}(r,x_i)I_{\varepsilon}(r,x_i)}{A_{\varepsilon}(r,x_i)}}\mds_{V_i}(.)dB_{x_i}(r)\Big\Vert_\mgamep^2 \Bigg) \nonumber \\
	&&\hspace{-4cm}+ 2C(\gamma)\E\Bigg( \Big\Vert \int_{\overline{\tau}}^{\overline{\tau}+\theta} \sum\limits_{\substack{i\, , \, j \\ x_i \sim x_j}} \dfrac{\sqrt{\mu_S S_{\varepsilon}(r,x_i)}}{\varepsilon}\Big( \mds_{V_j}(.)- \mds_{V_i}(.)\Big) dB_{x_i x_j}^S(r) \Big\Vert_\mgamep^2 \Bigg) \nonumber\\
	&&\hspace{-4cm} \le \dfrac{2C(\gamma)\overline{\beta}}{\varepsilon}\E\Bigg(\int_{\overline{\tau}}^{\overline{\tau}+\theta} \sum_m\sum_{i=1}^{\varepsilon^{-1}} \dfrac{S_{\varepsilon}(r,x_i)I_{\varepsilon}(r,x_i)}{A_{\varepsilon}(r,x_i)}\Big( \int_{V_i}\mathbf{f}_{m}^{\varepsilon}(x)\,dx \Big)^2 (1+\lambda_{m}^{\varepsilon})^{-\gamma} dr\Bigg) \nonumber \\
	&& \hspace{-4cm}+\dfrac{2C(\gamma)\mu_S}{\varepsilon}\E\Bigg(\int_{\overline{\tau}}^{\overline{\tau}+\theta} \sum_m \sum_{i=1}^{\varepsilon^{-1}}  S_{\varepsilon}(r, x_i) \Big( \int_{V_i} \nabla_{\varepsilon}^{\pm} \mathbf{f}_{m}^{\varepsilon}(x)\, dx \Big)^2(1+\lambda_{m}^{\varepsilon})^{-\gamma} dr\Bigg) \nonumber \\
	&& \hspace{-4cm}\le C(\overline{\beta} , \mu_S)\, \theta \longrightarrow 0, \qquad \text{as} \; \theta \to 0 .\nonumber
	\end{eqnarray}
	Hence  the condition \textbf{[A]} is proved.\\
	In the way, we prove similar estimates for $v_{\ep}$ and $w_{\ep}$. Then the process   $ \big\{ \Im_\ep(t),\; t\in[0,T], 0<\varepsilon<1 \big\}$ is  tight in   $C\big([0,T]; (\mrh^{-\gamma})^3\big)$, $\gamma>3/2$.
	\fpr
\end{pr}

\begin{lem}\label{convsub} For $3/2<\gamma <2$,
	the process $ \{ \, \Im_{\varepsilon}(t), \; t\in[0,T], \; 0<\varepsilon<1   \, \}$   converges in law in $C\big([0,T]\, ; \,  (\mrh^{-\gamma})^3\big)\cap  L^2\big(0,T; (\mrh^{-1})^3\big)$.
\end{lem}
\begin{pr} On the one hand, from \ref{tight},  the process $ \{ \, \Im_{\varepsilon}(t), \; t\in[0,T], \; 0<\varepsilon<1   \, \}$ is tight in $C\big([0,T]\, ; \, (\mrh^{-\gamma})^3\big) $, then along a subsequence, it converges in $C\big([0,T]\, ; \,  (\mrh^{-\gamma})^3\big) $. On the other hand the sequence $ \{ \, \Im_{\varepsilon}(t), \; t\in[0,T], \; 0<\varepsilon<1   \, \}$ is bounded in $ L^2\big(0,T\; ; \; (\mrh^{1-\gamma})^3\big)$. Indeed for all $\ep$ , we have 
	\begin{eqnarray}
	\mathbb{E}\bigg(\int_0^T\big\Vert u_{\ep}(t)\big\Vert_\ungam^2dt\bigg)&\le& C(\gamma)\mathbb{E}\bigg(\int_0^T\big\Vert u_{\ep}(t)\big\Vert_\ungamep^2dt\bigg)\; \; \quad  (\text{by using the inequality (\ref{equivl})})\n\\
	&=&C(\gamma)  \sum_m \mathbb{E}\Big(\int_0^T\langle\,u_{\ep}(t) , \f_{m}^{\ep}\,\rangle^2 dt\Big) (1+\lambda_{m,S}^{\ep})^{1-\gamma}\n \\
	&=& C(\gamma) \sum_m \mathbb{E}\Big(\int_0^T\langle\,u_{\ep}(t) , \f_{m}^{\ep}\,\rangle^2 dt\Big) (1+\lambda_{m,S}^{\ep})^{-\gamma}\n\\
	&+& C(\gamma) \sum_m \mathbb{E}\Big(\int_0^T\langle\,u_{\ep}(t) , \f_{m}^{\ep}\,\rangle^2 dt \Big)\lambda_{m,S}^{\ep}(1+\lambda_{m,S}^{\ep})^{-\gamma} \n \\
	&=& C(\gamma)\left\{\mathbb{E}\Big(\int_0^T \big\Vert u_{\ep}(t)\big\Vert_\mgamep^2 dt \Big)+ \mathbb{E}\Big(\int_0^T \big\Vert \nabla_{\!\!\ep}^+u_{\ep}(t)\big\Vert_\mgamep^2 dt \Big)\right\}\n\\
	&\le& C(\gamma)\left\{\mathbb{E}\Big(\int_0^T \big\Vert u_{\ep}(t)\big\Vert_\mgam^2 dt \Big)+ \mathbb{E}\Big(\int_0^T \big\Vert \nabla_{\!\!\ep}^+u_{\ep}(t)\big\Vert_\mgam^2 dt \Big)\right\}, \n
	\end{eqnarray}	
	where the third equality  follows from the fact that $$ \big\Vert \nabla_{\!\!\ep}^+u_{\ep}(t)\big\Vert_\mgamep^2= \sum_m \langle\,u_{\ep}(t) , \f_{m}^{\ep}\,\rangle^2 \lambda_{m,S}^{\ep}(1+\lambda_{m,S}^{\ep})^{-\gamma}\; \; (\text{see} \; \,\ref{an}(i)\; \, \text{in the Appendix below}).$$
	The inequality \eqref{nabl} ensures that $ \displaystyle \mathbb{E}\int_0^T \Big[\big\Vert u_{\ep}(t)\big\Vert_\mgam^2 dt +  \big\Vert \nabla_{\!\!\ep}^+u_{\ep}(t)\big\Vert_\mgam^2  \Big]dt$ is bounded by a constant independent of $\ep$. It then follows that 
	$$ \un{0<\ep<1}{\sup}\mathbb{E}\Big(\int_0^T\big\Vert u_{\ep}(t)\big\Vert_\ungam^2 dt\Big)\le C(\gamma).$$
	We have similar estimates for $ v_{\ep}$ and $w_{\ep}$. Thus
	\begin{eqnarray}
	\un{0<\ep<1}{\sup}\mathbb{E}\bigg(\int_0^T\big\Vert \Im_{\ep}(t)\big\Vert_\ungam^2 dt\bigg) &\le& C. \n
	\end{eqnarray}
	This implies that, from the sequence $ \{ \, \Im_{\varepsilon}(t), \; t\in[0,T], \; 0<\varepsilon<1   \, \}$, we can extract a subsequence which converges in law in $ L^2\big(0,T\; ; \;(\mrh^{1-\gamma})^3\big)$  endowed with the weak topology. 
	Furthermore, since the  imbedding of $\mrh^{1-\gamma}$ into $\mrh^{-1}$  is compact and we have the convergence in $C\left([0, T]\, ; (\mrh^{-\gamma})^3\right)$, then the extracted sequence converges in fact in $ L^2\big(0,T\; ; \; (\mrh^{-1})^3\big)$. Hence, we deduce that there exists a subsequence which converges in law in $C\big([0,T]\, ; \,  (\mrh^{-\gamma})^3\big)\cap  L^2\big(0,T\; ; \; (\mrh^{-1})^3\big)$.\\
	We note that the limit $ \displaystyle \Im :=(u, v, w)^{\mtt}$ of any convergent  subsequence satisfies the following system of stochastic PDEs
	\begin{equation}
	\left\{ 
	\begin{aligned} 
	d u(t) &= \mu_S\,\Delta u(t) dt+d\mathscr{M}^S(t)\\
	dv(t) &=\mu_I\,\Delta v dt+d\mathscr{M}^I(t) \\
	dw(t) &=\mu_R\,\Delta w(t)dt+d\mathscr{M}^R(t)
	\end{aligned}
	\right. 
	\end{equation}
	and  the solution of that system is unique. Then the whole process $ \{ \, \Im_{\varepsilon}(t), \; t\in[0,T], \; 0<\varepsilon<1   \, \}$ converges  in   $C\big([0,T]\, ; \,  (\mrh^{-\gamma})^3\big)\cap L^2\big(0,T\; ; \; (\mrh^{-1})^3\big)$.
	\fpr
\end{pr}

\begin{lem}\label{fu}
	As $\eps \to 0$, $f_\eps u_\eps \Longrightarrow fu$, and   $g_\eps v_\eps \Longrightarrow gv$  in $L^2\Big(0,T ; \mrh^{-1}\Big)$.
\end{lem}

\begin{pr}
	The convergence  $f_\eps u_\eps \Longrightarrow fu$ follows to the fact that $u_\eps \Longrightarrow u$ in $L^2\big(0,T\; ; \; \mrh^{-1}\big)$ and $f_\eps \longrightarrow f$ in $C\big([0,T]\; ; \; \mrh^{1}\big)$.  The proof of the convergence  $g_\eps v_\eps \Longrightarrow gv$   is similar. 
	\fpr
\end{pr}

%\begin{lem}
%$f_\ep u_\ep \longrightarrow fu$, \qquad  $g_\ep v_\ep \longrightarrow gv\; \; $ in $L^2(0, T; L^2(\tor))$.
%\end{lem}
%
%\begin{pr}
%Lemme de compacité de Thierry
%\fpr
%\end{pr}
We are now interested in the convergence of the process $\overline{\Im}_\ep:= \left(\overline{u}_\ep ,  \overline{v}_\ep, \overline{w}_\ep\right)$.
%{\color{blue} We have $f_\ep u_\ep \longrightarrow fu$, and  $g_\ep v_\ep \longrightarrow gv\; \; $ in $L^2(]0, T[\; ; \; L^2)$ ???}.
%We want to deduce that
%$(\overline{u}_\ep,\overline{v}_\ep,\overline{w}_\ep) \longrightarrow(\overline{u},\overline{v},\overline{w})$ in $(L^2(]0, T[\; ; \;L^2))^3$. 
\begin{lem}\label{estemps} For any $T>0$, there exists a positive constant $C$ such that
	\begin{eqnarray}
	\sup_{0\le t\le T}\Big(\big\Vert \overline{u}_\ep(t)\big\Vert_{L^2}^2+\big\Vert \overline{v}_\ep(t)\big\Vert_{L^2}^2 +\big\Vert \overline{w}_\ep(t)\big\Vert_{L^2}^2\Big)\n \\
	&& \hspace{-5cm} +C \int_{0}^{T}\Big(\big\Vert\nabla_{\varepsilon}^{+}\overline{u}_{\varepsilon}(s)\big\Vert_{L^2}^2 +\big\Vert\nabla_{\varepsilon}^{+}\overline{v}_{\varepsilon}(s)\big\Vert_{L^2}^2+\big\Vert\nabla_{\varepsilon}^{+}\overline{w}_{\varepsilon}(s)\big\Vert_{L^2}^2\Big) ds \le C\eta_{_T}e^{CT} , 
	\end{eqnarray}
	% and
	%\begin{eqnarray}\label{dt}
	%\int_0^T\Bigg(\Big\Vert \dfrac{d\,\overline{u}_\ep}{dt}(t)\Big\Vert^2_{\mrh^{-1}}+\Big\Vert \dfrac{d\,\overline{v}_\ep}{dt}(t)\Big\Vert^2_{\mrh^{-1}}+\Big\Vert \dfrac{d\,\overline{w}_\ep}{dt}(t)\Big\Vert^2_{\mrh^{-1}}\Bigg)dt
	%&\le &  C \eta_{_T}, 
	%\end{eqnarray}
	where $ \displaystyle\eta_{_T} :=\int_{0}^{T}\Big(\big\Vert u_\ep(s)\big\Vert_{\mrh^{-1}}^2 + \big\Vert v_\ep(s)\big\Vert_{\mrh^{-1}}^2\Big) ds.$
\end{lem}

\begin{pr} For all $t\in [0 , T],$  we have	
	\begin{eqnarray}
	\int_{0}^{t} \langle \dfrac{d\overline{u}_\ep}{ds}(s) \; , \; \overline{u}_\ep(s) \rangle ds \!\!\!&=&\!\!\!  \mu_S \int_{0}^{t} \langle\Delta_{\varepsilon} \overline{u}_{\varepsilon}(s)\; , \; \overline{u}_\ep(s) \rangle ds- \int_{0}^{t}\langle f_{\varepsilon}(s)\overline{u}_{\varepsilon}(s)\; , \; \overline{u}_\ep(s) \rangle ds\n \\ 
	&&\hspace{-2cm}-\int_{0}^{t}\langle g_{\varepsilon}(s)\overline{v}_{\varepsilon}(s)\; , \; \overline{u}_\ep(s) \rangle ds-\int_{0}^{t}\langle f_{\varepsilon}(s)u_{\varepsilon}(s)\; , \; \overline{u}_\ep(s) \rangle ds- \int_{0}^{t}\langle g_{\varepsilon}(s)v_{\varepsilon}(s)\; , \; \overline{u}_\ep(s) \rangle ds .\n
	\end{eqnarray}
	Then
	\begin{eqnarray}
	\big\Vert \overline{u}_\ep(t)\big\Vert_{L^2}^2 +2\mu_S \int_{0}^{t} \big\Vert\nabla_{\varepsilon}^{+}\overline{u}_{\varepsilon}(s)\big\Vert_{L^2}^2 ds &=& - 2\int_{0}^{t}\langle f_{\varepsilon}(s)\overline{u}_{\varepsilon}(s)\; , \; \overline{u}_\ep(s) \rangle ds- 2\int_{0}^{t}\langle g_{\varepsilon}(s)\overline{v}_{\varepsilon}(s)\; , \; \overline{u}_\ep(s) \rangle ds\n \\ 
	&-& 2\int_{0}^{t}\langle f_{\varepsilon}(s)u_{\varepsilon}(s)\; , \; \overline{u}_\ep(s) \rangle ds- 2\int_{0}^{t}\langle g_{\varepsilon}(s)v_{\varepsilon}(s)\; , \; \overline{u}_\ep(s) \rangle ds. \n
	\end{eqnarray}
	Since $\displaystyle f_\ep(t) u_\ep(t) \in \mrh^{-1}$  and  	$\displaystyle g_\ep(t) v_\ep(t) \in \mrh^{-1}$, then
	\begin{eqnarray}
	\big\Vert \overline{u}_\ep(t)\big\Vert_{L^2}^2 +2\mu_S \int_{0}^{t} \big\Vert\nabla_{\varepsilon}^{+}\overline{u}_{\varepsilon}(s)\big\Vert_{L^2}^2 ds &\le&  2\sup_{0\le s\le T} \Vert f_\ep(s)\Vert_{\infty} \int_{0}^{t}\big\Vert \overline{u}_\ep(s)\big\Vert_{L^2}^2 ds\n \\
	&& \hspace{-4cm}+\int_{0}^{t}\Big(\sup_{0\le s\le T} \Vert g_\ep(s)\Vert_{\infty}^2\big\Vert \overline{v}_\ep(s)\big\Vert_{L^2}^2+ \big\Vert \overline{u}_\ep(s)\big\Vert_{L^2}^2 \Big) ds\n \\ 
	&& \hspace{-4cm}+ 2\sup_{0\le s\le T}\big\Vert f_\ep(s)\big\Vert_{\mrh^{1,\ep}}\int_{0}^{t}\left[\big\Vert u_\ep(s)\big\Vert_{\mrh^{-1}}\big(\big\Vert \overline{u}_\ep(s)\big\Vert_{L^2}+ \big\Vert \nabla_{\!\!\ep}^+\overline{u}_\ep(s)\big\Vert_{L^2} \big) \right] ds\n \\
	&& \hspace{-4cm}+2\sup_{0\le s\le T}\big\Vert g_\ep(s)\big\Vert_{\mrh^{1,\ep}}\int_{0}^{t}\left[\big\Vert v_\ep(s)\big\Vert_{\mrh^{-1}}\big(\big\Vert \overline{u}_\ep(s)\big\Vert_{L^2}+ \big\Vert \nabla_{\!\!\ep}^+\overline{u}_\ep(s)\big\Vert_{L^2} \big) \right] ds.\n 
	\end{eqnarray}
	Let  $\delta$ be some constant such that $0<\delta <\dfrac{\mu_S}{C}.$ We have
	\begin{eqnarray}
	\big\Vert \overline{u}_\ep(t)\big\Vert_{L^2}^2 +2\mu_S \int_{0}^{t} \big\Vert\nabla_{\varepsilon}^{+}\overline{u}_{\varepsilon}(s)\big\Vert_{L^2}^2 ds &\le&  C\int_{0}^{t}\big\Vert \overline{u}_\ep(s)\big\Vert_{L^2}^2 ds
	+\int_{0}^{t}\Big(C\big\Vert \overline{v}_\ep(s)\big\Vert_{L^2}^2+ \big\Vert \overline{u}_\ep(s)\big\Vert_{L^2}^2 \Big) ds\n \\ 
	&&\hspace{-4cm}+ C\int_{0}^{t}\left[2\delta\big\Vert \overline{u}_\ep(s)\big\Vert_{L^2}^2+2\delta\big\Vert\nabla_{\!\!\ep}^+\overline{u}_\ep(s)\big\Vert_{L^2}^2+\dfrac{2}{\delta} \big\Vert u_\ep(s)\big\Vert_{\mrh^{-1}}^2 +\dfrac{2}{\delta} \big\Vert v_\ep(s)\big\Vert_{\mrh^{-1}}^2 \right] ds.\n 
	\end{eqnarray}
	Then
	\begin{eqnarray}\label{igu}
	\big\Vert \overline{u}_\ep(t)\big\Vert_{L^2}^2 +2(\mu_S-C\delta) \int_{0}^{t} \big\Vert\nabla_{\varepsilon}^{+}\overline{u}_{\varepsilon}(s)\big\Vert_{L^2}^2 ds &\le&  (C+2C\delta)\int_{0}^{t}\big\Vert \overline{u}_\ep(s)\big\Vert_{L^2}^2 ds\n 
	+C\int_{0}^{t}\big\Vert \overline{v}_\ep(s)\big\Vert_{L^2}^2\n \\
	&+&\dfrac{2C}{\delta}\int_{0}^{t}\Big( \big\Vert u_\ep(s)\big\Vert_{\mrh^{-1}}^2 + \big\Vert v_\ep(s)\big\Vert_{\mrh^{-1}}^2 \Big) ds.
	\end{eqnarray}
	In the same way, we prove that
	\begin{eqnarray}\label{igv}
	\big\Vert \overline{v}_\ep(t)\big\Vert_{L^2}^2 +2(\mu_I-C\delta) \int_{0}^{t} \big\Vert\nabla_{\varepsilon}^{+}\overline{v}_{\varepsilon}(s)\big\Vert_{L^2}^2 ds &\le&  (C+2C\delta)\int_{0}^{t}\big\Vert \overline{u}_\ep(s)\big\Vert_{L^2}^2 ds\n 
	+C\int_{0}^{t}\big\Vert \overline{v}_\ep(s)\big\Vert_{L^2}^2\n \\
	&+&\dfrac{2C}{\delta}\int_{0}^{t}\Big( \big\Vert u_\ep(s)\big\Vert_{\mrh^{-1}}^2 + \big\Vert v_\ep(s)\big\Vert_{\mrh^{-1}}^2 \Big) ds,
	\end{eqnarray}
	and 
	\begin{eqnarray}\label{igw}
	\big\Vert \overline{w}_\ep(t)\big\Vert_{L^2}^2 +2\mu_R\int_{0}^{t} \big\Vert\nabla_{\varepsilon}^{+}\overline{w}_{\varepsilon}(s)\big\Vert_{L^2}^2 ds &\le&  \dfrac{C}{\delta}\int_{0}^{t}\big\Vert v_\ep(s)\big\Vert_{\mrh^{-1}}^2 ds  \n \\
	&&\hspace{-1cm} +C\int_{0}^{t}\big\Vert \overline{w}_\ep(s)\big\Vert_{L^2}^2+C\int_{0}^{t}\big\Vert \overline{v}_\ep(s)\big\Vert_{L^2}^2 ds .
	\end{eqnarray}
	By adding the inequalities (\ref{igu}) , (\ref{igv}) and (\ref{igw}), we obtain
	\begin{eqnarray}
	\big\Vert \overline{u}_\ep(t)\big\Vert_{L^2}^2+\big\Vert \overline{v}_\ep(t)\big\Vert_{L^2}^2 +\big\Vert \overline{w}_\ep(t)\big\Vert_{L^2}^2+C \int_{0}^{t}\Big(\big\Vert\nabla_{\varepsilon}^{+}\overline{u}_{\varepsilon}(s)\big\Vert_{L^2}^2 +\big\Vert\nabla_{\varepsilon}^{+}\overline{v}_{\varepsilon}(s)\big\Vert_{L^2}^2+\big\Vert\nabla_{\varepsilon}^{+}\overline{w}_{\varepsilon}(s)\big\Vert_{L^2}^2\Big) ds \n \\
	&& \hspace{-12cm} \le C\int_{0}^{t}\Big(\big\Vert \overline{u}_\ep(s)\big\Vert_{L^2}^2+\big\Vert \overline{v}_\ep(s)\big\Vert_{L^2}^2+\big\Vert \overline{w}_\ep(s)\big\Vert_{L^2}^2 \Big)ds +C\int_{0}^{t}\Big(\big\Vert u_\ep(s)\big\Vert_{\mrh^{-1}}^2 + \big\Vert v_\ep(s)\big\Vert_{\mrh^{-1}}^2\Big) ds. \n
	\end{eqnarray}
	Hence  applying Gronwall's Lemma we obtain
	\begin{eqnarray}
	\sup_{0\le t\le T}\Big(\big\Vert \overline{u}_\ep(t)\big\Vert_{L^2}^2+\big\Vert \overline{v}_\ep(t)\big\Vert_{L^2}^2 +\big\Vert \overline{w}_\ep(t)\big\Vert_{L^2}^2\Big)\n \\
	&& \hspace{-5cm} +C \int_{0}^{T}\Big(\big\Vert\nabla_{\varepsilon}^{+}\overline{u}_{\varepsilon}(s)\big\Vert_{L^2}^2 +\big\Vert\nabla_{\varepsilon}^{+}\overline{v}_{\varepsilon}(s)\big\Vert_{L^2}^2+\big\Vert\nabla_{\varepsilon}^{+}\overline{w}_{\varepsilon}(s)\big\Vert_{L^2}^2\Big) ds \le C\eta_Te^{CT}.  \n
	\end{eqnarray}
	\hfill $\square$
\end{pr}

We want to deduce from the fact that the pair $(u_\eps,v_\eps)$ converges in law towards $(u,v)$ in $L^2(0,T; (\mrh^{-1})^2)$,
the convergence in law of $(\overline{u}_\eps,\overline{v}_\eps,\overline{w}_\eps)$.

\begin{lem}\label{convbar} The process $\{\, \left(\overline{u}_\ep(t),\overline{v}_\ep(t), \overline{w}_\ep(t)\right), 0\le t\le T, \; \; 0<\eps<1\, \}\Rightarrow \{\, (\overline{u}(t),\overline{v}(t), \overline{w}(t)), 0\le t\le T\, \}$ in 
	$L^2\big(0,T; (L^2)^3\big)\cap C([0,T];(\mrh^{-1})^3)$, where the limit $\{\,\left(\overline{u}(t),\overline{v}(t), \overline{w}(t)\right), 0\le t\le T\,\}$ is the unique solution of  the following  system of parabolic PDEs
	\begin{equation}\label{bar_EDP}
	\left\{ 
	\begin{aligned} 
	\dfrac{d\overline{u}}{dt}(t) &=  \mu_S\,\Delta \overline{u} (t)-f(t)\overline{u}(t)-g(t)\overline{v}(t)-f(t)u(t)-g(t)v(t)\\
	\dfrac{d\overline{v}}{dt}(t) &=  \mu_I\,\Delta \overline{v}(t)+ f(t)\overline{u}(t)+g(t)\overline{v}(t)+f(t)u(t)+g(t)v(t)-\alpha\left(v(t)+\overline{v}(t\right)
	\\
	\dfrac{d\overline{w}}{dt}(t) &=  \mu_R\,\Delta \overline{w}(t)+ \alpha\left(v(t)+\overline{v}(t)\right)
	\\ 
	\overline{u}(0)&=\overline{v}(0)=\overline{w}(0)=0.
	\end{aligned}
	\right. 
	\end{equation}
\end{lem}

\begin{pr}
	Let 
	\begin{align*}
	\overline{\Im}_\eps(t)&=\begin{pmatrix}\overline{u}_\eps(t)\\ \overline{v}_\eps(t)\\ \overline{w}_\eps(t)\end{pmatrix},\quad
	F_\eps(t)=\begin{pmatrix} -f_\eps(t)u_\eps(t)-g_\eps(t)v_\eps(t)\\
	f_\eps(t)u_\eps(t)+(g_\eps(t)-\alpha) v_\eps(t)\\ \alpha v_\eps(t)\end{pmatrix},\\[4mm]
	\Lambda_\eps(t)&=  \begin{pmatrix} \mu_S\Delta_\eps-f_\eps(t) & -g_\eps(t) & 0\\
	f_\eps(t) & \mu_I\Delta_\eps+g_\eps(t)-\alpha& 0\\
	0  & 0 & \mu_R\Delta_\eps+\alpha \end{pmatrix}.
	\end{align*}
	Note that both $\overline{\Im}_\eps$ and $F_\eps$ belong to $L^2(0,T ; (\hd)^3)$.                
	We have the following system of ODEs
	\begin{align}\label{EDP_eps}
	\frac{d\overline{\Im}_\eps}{dt}(t)=\Lambda_\eps(t)\overline{\Im}_\eps(t)+F_\eps(t),\quad \overline{\Im}_\eps(0)=0\,.
	\end{align}                                
	\ref{fu}  tells us that whenever, as $\eps \to 0$, 
	\begin{align*}
	F_\eps \Longrightarrow F \ \text{ in } L^2(0,T ; (\mrh^{-1})^3),
	\end{align*}
	where 
	\begin{align}
	F(t)=\begin{pmatrix} -f(t)u(t)-g(t)v(t)\\
	f(t)u(t)+g(t)v(t)-\alpha v(t)\\ \alpha v(t)\end{pmatrix}\,.
	\end{align}  
	We apply the well--known theorem due to Skorohod, which asserts that redefining the probability space, we can assume that $F_\eps\to F$ a.s. strongly in
	$L^2((0,T); (\mrh^{-1})^3)$.  Our assumptions and the hypotheses
	imply that both $\overline{\Im}_\eps$ and $\nabla_\eps^+\overline{\Im}_\eps$ are bounded
	in $(L^2((0,T)\times \tor))^3$. Hence along a subsequence $\overline{\Im}_\eps\to\overline{\Im}$ and 
	$\displaystyle\nabla_\eps^+\overline{\Im}_\eps\to\overline{G}$ in $(L^2((0,T)\times \tor))^3$ weakly. However, it follows from a duality argument that  
	$ \overline{G}=\nabla\overline{\Im}$, and taking the weak limit in  \eqref{EDP_eps},
	we deduce that $\overline{\Im}$ is the unique solution of  the system of parabolic PDEs 
	\begin{align*}
	\frac{d \overline{\Im}}{dt}(t)=\Lambda(t)\overline{\Im}(t)+F(t),\quad \overline{\Im}(0)=0\,,   
	\end{align*}                                       
	with
	\begin{align}
	\Lambda(t)= \begin{pmatrix} \mu_S\Delta-f(t) & -g(t) &0 \\
	f(t) & \mu_I\Delta+g(t)-\alpha & 0 \\
	0&0 & \mu_R\Delta+\alpha \end{pmatrix} .
	\end{align} 
	Hence all converging subsequences have the same limit, and the whole sequence converges.  
	
	We now show that the pair $\langle\overline{\Im}_\eps,\nabla_\eps^+\overline{\Im}_\eps\rangle$ converges strongly in 
	$(L^2((0,T)\times \T^1))^6$. We first note that  both $\overline{\Im}_\eps$ and 
	$\nabla_\eps^+\overline{\Im}_\eps$ are bounded in $(L^2((0,T)\times \T^1))^3$, but also 
	$\frac{d}{dt}\overline{\Im}_\eps$ is bounded in $L^2((0,T);(\mrh^{-1}(\T^1))^3)$. From these estimates, we deduce
	with the help of Theorem 5.4 in Droniou et al. (2018) that $\overline{\Im}_\eps\to\overline{\Im}$ strongly in $(L^2((0,T)\times \T^1))^3$. Next we deduce from \eqref{EDP_eps} that 
	\[ \frac{1}{2}\frac{d\big\Vert\overline{\Im}_\eps(t)\big\Vert_{L^2}^2}{dt}=\langle\Lambda_\eps\overline{\Im}_\eps(t),\overline{\Im}_\eps(t)\rangle+\langle F_\eps(t),\overline{\Im}_\eps(t)\rangle,\] hence 
	\begin{align}
	&\frac{1}{2}\big\Vert\overline{\Im}_\eps(T)\big\Vert_{L^2}^2+\int_0^T\left[\mu_S\big\Vert\nabla_\eps^+\overline{u}_\eps(t)\big\Vert_{L^2}^2+\mu_I\big\Vert\nabla_\eps^+\overline{v}_\eps(t)\big\Vert_{L^2}^2+\mu_R\big\Vert\nabla_\eps^+\overline{w}_\eps(t)\big\Vert_{L^2}^2\right]dt\label{ident_eps}\\
	&=\int_0^T\Big[\langle f_\eps(t)\overline{u}_\eps(t)+g_\eps(t)\overline{v}_\eps(t),\overline{v}_\eps(t)-\overline{u}_\eps(t)\rangle+\big\Vert\sqrt{\alpha}\;\overline{w}_\eps(t)\big\Vert_{L^2}^2-\big\Vert\sqrt{\alpha}\;\overline{v}_\eps(t)\big\Vert_{L^2}^2+\langle F_\eps(t),\overline{\Im}_\eps(t)\rangle\Big]dt\,.\nonumber
	\end{align}
	We have an analogous identity for the limiting quantities, namely:
	\begin{align}
	&\frac{1}{2}\big\Vert\overline{\Im}(T)\big\Vert_{L^2}^2+\int_0^T\left[\mu_S\big\Vert\nabla\overline{u}(t)\big\Vert_{L^2}^2+\mu_I\big\Vert\nabla\overline{v}(t)\big\Vert_{L^2}^2+\mu_R\big\Vert\nabla\overline{w}(t)\big\Vert_{L^2}^2\right]dt\label{ident}\\
	&=\int_0^T\Big[\langle f(t)\overline{u}(t)+g(t)\overline{v}(t),\overline{v}(t)-\overline{u}(t)\rangle+\big\Vert\sqrt{\alpha}\;\overline{w}(t)\big\Vert_{L^2}^2-\big\Vert\sqrt{\alpha}\;\overline{v}(t)\big\Vert_{L^2}^2+\langle F(t),\overline{\Im}(t)\rangle\Big]dt\,.\nonumber
	\end{align}
	It follows from the strong convergence of $F_\eps$ to $F$ in $L^2(0,T;(\mrh^{-1})^3)$, the strong convergence of $\overline{\Im}_\eps\to\overline{\Im}$ in $(L^2((0,T)\times \T^1))^3$ and the weak convergence of $\nabla_\eps^+\overline{\Im}_\eps$ to $\nabla\overline{\Im}$ in $(L^2((0,T)\times \T^1))^3$ that the right hand side of \eqref{ident_eps} converges to the right hand side of \eqref{ident}. Hence the left hand side of \eqref{ident_eps} converges to the left hand side of \eqref{ident}. Consequently
	\begin{eqnarray}\label{strong_conv}
	\frac{1}{2}\big\Vert\overline{\Im}_\eps(T)-\overline{\Im}(T)\big\Vert_{L^2}^2&+&\int_0^T\left[\mu_S\big\Vert\nabla_\eps^+\overline{u}_\eps(t)-\nabla\overline{u}(t)\big\Vert_{L^2}^2+\mu_I\big\Vert\nabla_\eps^+\overline{v}_\eps(t)-\nabla\overline{v}(t)\big\Vert_{L^2}^2 \right.\nonumber \\
	&&\hspace{2.5cm}+\left.\mu_R\big\Vert\nabla_\eps^+\overline{w}_\eps(t)-\nabla\overline{w}(t)\big\Vert_{L^2}^2\right]dt\to 0\,.
	\end{eqnarray}
	This last result follows from the convergence of the left hand side of \eqref{ident_eps} to that of \eqref{ident}, and
	the facts that 
	\begin{eqnarray}
	\langle\overline{\Im}_\eps(T),\overline{\Im}(T)\rangle&\to &\big\Vert\overline{\Im}(T)\big\Vert_{L^2}^2,\n
	\end{eqnarray}
	and
	\begin{eqnarray}
	\int_0^T\left[\mu_S\langle\nabla_\eps^+\overline{u}_\eps(t),\nabla\overline{u}(t)\rangle+\mu_I\langle\nabla_\eps^+\overline{v}_\eps(t),\nabla\overline{v}(t)\rangle 
	+\mu_R\langle\nabla_\eps^+\overline{w}_\eps(t),\nabla\overline{w}(t)\rangle\right]dt&&\n\\
	&& \hspace{-6cm} \to
	\int_0^T\left[\mu_S\big\Vert\nabla\overline{u}(t)\big\Vert_{L^2}^2+\mu_I\big\Vert\nabla\overline{v}(t)\big\Vert_{L^2}^2+\mu_R\big\Vert\nabla\overline{w}(t)\big\Vert_{L^2}^2\right]dt\,. \n
	\end{eqnarray}
	The second convergence follows from the fact that $\nabla_\eps^+\overline{\Im}_\eps\to\nabla\overline{\Im}$ in $(L^2((0,T)\times \T^1))^3$ weakly. Concerning the first one, we deduce from the equations and the above statements that $\overline{\Im}_\eps(T)\to\overline{\Im}(T)$ weakly in $(\mrh^{-1})^3$. But since that sequence is bounded in $(L^2(\T^1))^3$, it also converges weakly in $(L^2(\T^1))^3$.
	
	The fact that $\nabla_\eps^+\overline{\Im}_\eps\to\nabla\overline{\Im}$ strongly in $(L^2((0,T)\times \T^1))^3$ clearly follows from
	\eqref{strong_conv}.
	
	The above arguments imply that a.s.
	\[ \langle\overline{\Im}_\eps,\nabla_\eps^+\overline{\Im}_\eps\rangle\to \langle\overline{\Im},\nabla\overline{\Im}\rangle \ \text{ strongly in }(L^2((0,T)\times \T^1))^6\,.\]
	Now the convergence $\overline{\Im}_\eps\to\overline{\Im}$ in $C([0,T];(H^{-1})^3)$ follows readily from the equation. 
	\fpr
\end{pr}

\ref{convsub} says that $\Im_\ep\Rightarrow \Im$ in $C\big([0,T]\, ; \,  (\mrh^{-\gamma})^3\big)\cap  L^2\big(0,T; (\mrh^{-1})^3\big)$, we have used in \ref{convbar} the Skorohod theorem to deduce that $\overline{\Im}_\ep\Rightarrow
\overline{\Im}$ in $L^2(0,T(L^2)^3)\cap C([0,T];(\mrh^{-\gamma})^{-1})$. Hence the same Skorohod theorem allows us to take the limit in the sum $\Im_\ep+\overline{\Im}_\ep$, which yields the following result.
\begin{thm}[Functional central limit theorem]\label{tcl}
	For $3/2<\gamma <2$, as $\ep \to 0$, $\{\mathscr{Y}_{\ep}(t)\; , \;  0\le t\le T \}_{0<\ep<1} \Longrightarrow \{\mathscr{Y}(t)\; , \;  0\le t\le T \}$ in $C\big([0,T]\, ; \,  (\mrh^{-\gamma})^3\big)\cap  L^2\big(0,T; (\mrh^{-1})^3\big)$, where the limit $\mathscr{Y}$ is solution of the following system of SPDEs : for all $\varphi \in \mrh^{1}$
	\begin{equation}
	\left\{
	\begin{aligned}
	\big\langle \,  \mcru(t), \varphi \ \big\rangle_{\!{\mrh^{-1},\mrh^{1}}}
	&=\mu_S\int_0^t\!\! \big\langle \;  \mcru(r)\; ,\; \Delta\varphi \; \big\rangle_{\!{\mrh^{-1},\mrh^{1}}} dr +\!\! \int_0^t\!\! \big\langle \, \mcrv(r) , \,\beta(.) \dfrac{\mathbf{i}(r)\big(\mathbf{i}(r)+\mathbf{r}(r)\big)}{\mathbf{a}^2(r)}\varphi \, \big\rangle_{\!{\mrh^{-1},\mrh^{1}}} dr\\
	&\hspace{-0.5cm}+ \int_0^t \big\langle \; \mcru(r), \, \beta(.) \dfrac{\mathbf{s}(r)\big(\mathbf{s}(r)+\mathbf{r}(r)\big)}{\mathbf{a}^2(r)}\varphi \; \big\rangle_{\!{\mrh^{-1},\mrh^{1}}} dr + \big\langle\, \scm^S(t)\, , \, \varphi\; \big\rangle_{\!{\mrh^{-1},\mrh^{1}}} 
	\\[2mm]
	\big\langle \;  \mcrv(t), \varphi \; \big\rangle_{\!{\mrh^{-1},\mrh^{1}}} &= \mu_I\int_0^t \!\!\big\langle \,  \mcrv(r)\; ,\; \Delta\varphi \, \big\rangle_{\!{\mrh^{-1},\mrh^{1}}} dr - \int_0^t\big\langle \; \mcrv(r), \, \beta(.) \dfrac{\mathbf{i}(r)\big(\mathbf{i}(r)+\mathbf{r}(r)\big)}{\mathbf{a}^2(r)}\varphi \; \big\rangle_{\!{\mrh^{-1},\mrh^{1}}}  dr\\
	& \hspace{-2cm}- \int_0^t \big\langle \;  \mcru(r), \, \beta(.) \dfrac{\mathbf{s}(r)\big(\mathbf{s}(r)+\mathbf{r}(r)\big)}{\mathbf{a}^2(r)}\varphi \; \big\rangle_{\!{\mrh^{-1},\mrh^{1}}} dr + \int_0^t\big\langle \; \mcrv(r), \alpha(.)\varphi \; \big\rangle_{\!{\mrh^{-1},\mrh^{1}}} dr+\big\langle\, \scm^I(t)\, , \, \varphi\; \big\rangle_{\!{\mrh^{-1},\mrh^{1}}} 
	\\[3mm] 
	\big\langle \;  \mcrw(t), \varphi \; \big\rangle_{\!{\mrh^{-1},\mrh^{1}}}
	&= \mu_R \int_0^t \big\langle \;  \mcrw(r)\,,\,\Delta\varphi \; \big\rangle_{\!{\mrh^{-1},\mrh^{1}}} dr - \int_0^t\big\langle \; \mcrv(r), \alpha(.)\varphi \; \big\rangle_{\!{\mrh^{-1},\mrh^{1}}} dr + \big\langle\; \scm^R(t)\, , \, \varphi\; \big\rangle_{\!{\mrh^{-1},\mrh^{1}}} \, .
	\end{aligned}
	\right.
	\end{equation}
	
\end{thm}

\medspace

\textbf{\underline{Final remarks}:}  $\bullet$ Our functional central limit theorem is established in dimension 1. The difficulty in higher dimension is  the following. $\gamma>3/2$ has to be replaced by $\gamma>1+d/2$.  Then in \ref{convsub} we have convergence in $L^2(0,T;(H^{1-\gamma})^3)\cap C([0,T];(H^{-\gamma})^3)$. Note that $1-\gamma<-d/2$. Already in dimension 2, we have $1-\gamma<-1$, and there is a serious difficulty with the analog of \ref{convbar}.
\\
$\large \bullet$ In this work, we have first let $\mathbf{N}\to\infty$, while $\eps>0$ is fixed, and then let $\eps\to0$.
The case where $\mathbf{N}\to +\infty$ and $\eps \to 0$ together, with some constraint on the relative speeds of convergence (which does not allow $\mathbf{N}$ to converge too slowly to $\infty$ while $\eps\to0$)  will be the subject of future work.

\bigskip

\bmhead{Acknowledgments}

The author T\'enan Yeo would like to thank the Marseille Mathematics Institute (I2M)
for funding his stay in Marseille, during which part of this work was carried out.

\section*{Declarations}

\bmhead{Funding} The authors were supported by  their respective university.

\bmhead{Conflict of interest/Competing interests} The authors declare that they have no conflict of interest and no competing interests.

The authors have no competing interests to declare that are relevant to the content of this article.

\begin{appendices}
	
\section{}

\begin{lem}
	Let $(h_\ep)_{0<\ep<1}$ be a sequence of $\mathtt{H}_{\varepsilon}$. If $(h_\ep)_{0<\ep<1}$ is bounded in $\mrh^{1,\ep}$, then it is relatively compact in $L^2$, and the limit of any convergent subsequence belongs to $\mrh^{1}$.
\end{lem}	
\begin{pr} By using the fact that  the sequence $(h_\ep)$ is bounded in $L^2$ and $\big\Vert \nabla_{\!\!\ep}^+h_\ep\big\Vert_{L^2}\le C\big\Vert h_\ep \big\Vert_{_{\mrh^{1}}}$,  then the result of the compactness follows from the compactness theorem of Kolmogorov in $L^2$.
	
	The fact the limit of  any convergent subsequence belong to $\mrh^{1}$, follows from the discrete integrating by part
	$$\int_{\tor}\nabla_{\!\!\ep}^+h_\ep(x)\varphi(x)dx= -\langle\, \int_0^. h_\ep(y)dy\; , \; \nabla_{\!\!\ep}^+\varphi \; \rangle, $$ and letting $\ep$ go to zero in this equation. \fpr
\end{pr}

\begin{lem}\label{an} ~\
	%\begin{enumerate}[label=\bf (\roman*)]
	For all $\displaystyle u_{\ep} \in  \mathtt{H}_{\varepsilon} $
	
	$ \displaystyle \big\Vert \nabla_{\!\!\ep}^- u_{\ep}\big\Vert_\mgamep^2 = \big\Vert \nabla_{\!\!\ep}^+ u_{\ep}\big\Vert_\mgamep^2= \sum_m \left( \langle\,u_{\ep} , \varphi_{m}^{\ep}\,\rangle^2 + \langle\,u_{\ep} , \psi_{m}^{\ep}\,\rangle^2   \right) \lambda_{m}^{\ep}(1+\lambda_{m}^{\ep})^{-\gamma}. $
	
	%For all $\displaystyle u \in \mrh^{1-\gamma}$
	
	%$\displaystyle \big\Vert \nabla u \big\Vert_{_{\!{\mrh^{-\gamma}}}}^2+\Vert u \Vert_{_{\!{\mrh^{-\gamma}}}}^2
	%	=\Vert u \Vert_{_{\!{\mrh^{1-\gamma}}}}^2.$
	%	\end{enumerate}
\end{lem}

\begin{pr} We have
	$$ 	\nabla^-_{\!\!\ep} \varphi_m^{\ep} = -b_{m,\ep} \varphi_m^{\ep}- a_{m,\ep}\psi_m^{\ep} \;  \; \text{ and } \; \; 
	\nabla^-_{\!\!\ep} \psi_m^{\ep} = a_{m,\ep} \varphi_m^{\ep}-b_{m,\ep}\psi_m^{\ep}, $$
	where 	 $a_{m,\ep}=\ep^{-1}\sin(\pi m\ep)$ and  $b_{m,\ep}=\ep^{-1} (\cos(\pi m\ep)-1)$.
	
	We have $$ a_{m,\ep}^2+b_{m,\ep}^2=\lambda_m^{\ep}.$$
	
	Let $\displaystyle u_{\ep} \in  \mathtt{H}_{\varepsilon} $. We have
	%\begin{eqnarray}
	%\big\Vert \nabla_{\!\!\ep}^- u_{\ep}\big\Vert_{-\gamma,\ep}^2&=&\sum_m\left( \langle\, u_{\ep},\nabla_{\!\!\ep}^+ \varphi_m\, \rangle^2+ \langle\, u_{\ep} , \nabla_{\!\!\ep}^+ \psi_m \,\rangle^2\right)(1+\lambda_m^{\ep})^{-\gamma}\n \\
	%&=&\sum_m\left(\langle\, u_{\ep},b_{m,\ep}\varphi_m^{\ep}-a_{m,\ep} \psi_m^{\ep} \,\rangle^2+\langle\, u_{\ep}, a_{m,\ep} \varphi_m^{\ep}+b_{m,\ep} \psi_m^{\ep} \,\rangle^2\right)(1+\lambda_m^{\ep})^{-\gamma} \n \\
	%&=&\sum_m\left([b_{m,\ep} \langle\, u_{\ep},\varphi_m^{\ep}\,\rangle -a_{m,\ep}\langle\, u_{\ep}, \psi_m^{\ep} \,\rangle ]^2+[a_{m,\ep}\langle\, u_{\ep},\varphi_m^{\ep}\,\rangle+b_{m,\ep} \langle\, u_{\ep}, \psi_m^{\ep}\,\rangle ]^2\right)(1+\lambda_m^{\ep})^{-\gamma} \n\\
	%&=&\sum_m\left([a_{m,\ep}^2+b_{m,\ep}^2]\{\langle\, u_{\ep},\varphi_m^{\ep}\,\rangle^2+\langle\, u_{\ep}, \psi_m^{\ep}\,\rangle^2\}\right)(1+\lambda_m^{\ep})^{-\gamma}\n \\
	%&=&\sum_m\left(\langle\, u_{\ep}, \varphi_m^{\ep}\,\rangle^2+\langle\, u_{\ep}, \psi_m^{\ep}\,\rangle^2\right)\lambda_m^\ep(1+\lambda_m^{\ep})^{-\gamma} .\n
	%\end{eqnarray}	
	\begin{eqnarray}
	\big\Vert \nabla_{\!\!\ep}^+ u_{\ep}\big\Vert_\mgamep^2&=&\sum_m\left( \langle\, u_{\ep},\nabla_{\!\!\ep}^- \varphi_m\, \rangle^2 + \langle\, u_{\ep} , \nabla_{\!\!\ep}^- \psi_m \,\rangle^2\right)(1+\lambda_m^{\ep})^{-\gamma}\n \\
	&=&\sum_m\left(\langle\, u_{\ep}, -b_{m,\ep}\varphi_m^{\ep}-a_{m,\ep} \psi_m^{\ep} \,\rangle^2+\langle\, u_{\ep}, a_{m,\ep} \varphi_m^{\ep}-b_{m,\ep} \psi_m^{\ep} \,\rangle^2\right)(1+\lambda_m^{\ep})^{-\gamma} \n \\
	&=&\sum_m\left([-b_{m,\ep} \langle\, u_{\ep},\varphi_m^{\ep}\,\rangle -a_{m,\ep}\langle\, u_{\ep}, \psi_m^{\ep} \,\rangle ]^2+[a_{m,\ep}\langle\, u_{\ep},\varphi_m^{\ep}\,\rangle-b_{m,\ep} \langle\, u_{\ep}, \psi_m^{\ep}\,\rangle ]^2\right)(1+\lambda_m^{\ep})^{-\gamma} \n\\
	&=&\sum_m\left([a_{m,\ep}^2+b_{m,\ep}^2]\{\langle\, u_{\ep},\varphi_m^{\ep}\,\rangle^2+\langle\, u_{\ep}, \psi_m^{\ep}\,\rangle^2\}\right)(1+\lambda_m^{\ep})^{-\gamma}\n \\
	&=&\sum_m\left(\langle\, u_{\ep}, \varphi_m^{\ep}\,\rangle^2+\langle\, u_{\ep}, \psi_m^{\ep}\,\rangle^2\right)\lambda_m^\ep(1+\lambda_m^{\ep})^{-\gamma} .\n
	\end{eqnarray}
	The proof of   
	$\displaystyle \big\Vert \nabla_{\!\!\ep}^- u_{\ep}\big\Vert_\mgamep^2 = \sum_m \left( \langle\,u_{\ep} , \varphi_{m}^{\ep}\,\rangle^2 + \langle\,u_{\ep} , \psi_{m}^{\ep}\,\rangle^2   \right) \lambda_{m}^{\ep}(1+\lambda_{m}^{\ep})^{-\gamma} $
	is similar by noting that
	$$\nabla^+_{\!\!\ep} \varphi_m^{\ep} = b_{m,\ep} \varphi_m^{\ep}- a_{m,\ep}\psi_m^{\ep} \;  \; \text{ and } \; \; 
	\nabla^+_{\!\!\ep} \psi_m^{\ep} = a_{m,\ep} \varphi_m^{\ep}+b_{m,\ep}\psi_m^{\ep}.$$
	%		Let $u \in  \mrh^{1-\gamma}$. 	
	%		\begin{eqnarray}
	%		\big\Vert \nabla u \big\Vert_{_{\!{\mrh^{-\gamma}}}}^2&=&\sum_{m} \Big(\langle  u, \nabla\varphi_{m} \rangle^2+ \langle  u, \nabla\psi_{m} \rangle^2\Big) (1+\lambda_{m})^{-\gamma} \n \\
	%		&=&\sum_{m} \Big( \langle  u, \varphi_{m} \rangle^2+ \langle  u, \psi_{m} \rangle^2\Big)\pi^2m^2(1+\lambda_{m})^{-\gamma} \n \\
	%		&=& \sum_{m} \langle  u, \f_{m} \rangle^2\pi^2 m^2(1+\lambda_{m})^{-\gamma}.\n 
	%		\end{eqnarray}
	%		Then
	%		\begin{eqnarray}
	%		\big\Vert \nabla u \big\Vert_{_{\!{\mrh^{-\gamma}}}}^2+\Vert u \Vert_{_{\!{\mrh^{-\gamma}}}}^2
	%		&=& \sum_{m} \langle  u, \f_{m} \rangle^2\pi^2 m^2(1+\lambda_{m})^{-\gamma}+\sum_m \langle  u, \f_{m} \rangle^2(1+\lambda_{m})^{-\gamma}\n \\
	%		&=& \sum_{m}  \langle  u, \f_{m} \rangle^2(1+\lambda_{m})^{1-\gamma}\n\\
	%		&=&\Vert u \Vert_{_{\!{\mrh^{1-\gamma}}}}^2. \n
	%		\end{eqnarray}
	
	\fprb
\end{pr}





\end{appendices}



\begin{thebibliography}{2}
	%\rhead{References}
	\bibitem{hj} H. Bahouri, J. Chemin, R. Danchin. \textit{Fourier Analysis and Nonlinear Partial
		Differential Equations}, Springer (2011).
	\bibitem{db93} D. J. Blount.  Limit theorems for a sequence of nonlinear reaction-diffusion systems.\textit{ Stochastic Processes and their Applications}, \textbf{45}(2), 193-207 (1993).
	%\bibitem{db87} D. J. Blount. Comparison of a stochastic model of a chemical reaction with diffusion and the deterministic model. Ph.d., The University of Wisconsin-Madison, 1987.
	\bibitem{BP2019} T. Britton, E. Pardoux.  \textit{Stochastic Epidemic Models with Inference}, Lecture
	notes in Math. \textbf{2255}, Springer (2019).
	\bibitem{th} J. Droniou, R. Eymard, T.  Gallou\"et, R. Herbin. \textit{The gradient discretisation method}. Springer (2018).
	\bibitem{jm86} A. Joffe, M. Metivier.  Weak Convergence of Sequences of Semimartingales with Applications to Multitype Branching Processes, \textit{Advances in Applied Probability}, Vol. \textbf{18}, No. 1 (Mar., 1986), pp. 20-65.
	\bibitem{mc} W.O. Kermack, A.G. McKendrick. Proc. Roy. Soc. A 115
	, 700. Reprinted in \textit{ Bull. Math. Biol.} \textbf{53} (1991) 33.
	\bibitem{kg} P. Kotelenez.  Gaussian approximation to the nonlinear reaction-diffusion equation.
	Report 146, Universit\"{a}t Bremen Forschungsschwerpunkt Dynamische Systemes (1986).
	\bibitem{kot84} P. Kotelenez. A stopped Doob inequality for stochastic convolution integrals and stochastic evolution equations, Stochastic analysis and applications, \textbf{2}(3), 245-265 (1984).
	\bibitem{tg} T. G. Kurtz. Limit theorems for sequences of of jump Markov processes approximating ordinary differential processes. J. \textit{Appl. Probab.}, \textbf{8}: 344-356 (1971).
	%\bibitem{Kur81}T.G. Kurtz. Approximation of population processes. CBMS-NSF Regional Conference Series in Applied Mathematics, 36 SIAM, Philadelphia (1981).
	\bibitem{PY2018} M. N'zi, E. Pardoux, T. Yeo. A SIR model on a refining spatial grid I. Law of large numbers. \textit{Applied Mathematics $\&$ Optimization} 83 : 1153-1189 (2021).
	%\bibitem{p75} E. Pardoux. Equations aux d\'eriv\'ees partielles stochastiques non lin\'eaires monotones; Etude de solutions fortes de type It\^o. Th\'ese, Univ. Paris Sud, 1975.
	\bibitem{tay} M. E. Taylor. \textit{Partial Differential Equations III Nonlinear Equations}. Springer (1991).
	
	
	
	%\bibitem{pazy} Pazy, A.  \textit{Semigroups of linear operators and applications to partial differential equations}, (Vol. 44). Springer Science $\&$ Business Media (2012).
	%\bibitem{b87} D. J. Blount. Comparison of a stochastic model of a chemical reaction with diffusion
	%and the deterministic model. Ph.d., The University of Wisconsin-Madison,
	%1987.
\end{thebibliography}
\end{document}